\documentclass[article,hidelinks,onefignum,onetabnum]{siamart220329}

\usepackage{algorithmic}
\usepackage{amsmath}
\usepackage{graphicx}
\usepackage{booktabs}
\usepackage{soul}
\usepackage{subcaption}

\usepackage{todonotes}
\makeatletter 
\@mparswitchfalse%
\makeatother
\normalmarginpar 

\usepackage{tikz}
\usetikzlibrary{arrows}
\usetikzlibrary{decorations.pathreplacing}
\usetikzlibrary{shapes}

\renewcommand{\vec}[1]{\boldsymbol #1}
\newcommand{\vv}{\vec{v}}
\newcommand{\vu}{\vec{u}}
\newcommand{\ve}{\vec{e}}
\newcommand{\vr}{\vec{r}}
\newcommand{\vf}{\vec{f}}
\newcommand{\vg}{\vec{g}}
\newcommand{\vtau}{\vec{\tau}}
\newcommand{\fg}{\Omega^{h}}
\newcommand{\cg}{\Omega^{mh}}

\newcommand{\TheTitle}{
Multigrid Reduction in Time for Chaotic Dynamical Systems
}
\newcommand{\TheAuthors}{
  D. A. Vargas,
  R. D. Falgout,
  S. G{\"u}nther,
  J. B. Schroder
}

\newcommand{\TheFunding}{%
This work was performed under the auspices of the US Department of
Energy by Lawrence Livermore National Laboratory under Contract
DE-AC52-07NA27344 (LLNL-JRNL-838414).
}

\ifpdf{}
\hypersetup{
  pdftitle={\TheTitle},
  pdfauthor={\TheAuthors}
}
\fi

\title{{\TheTitle}\thanks{\TheFunding}}
\author{\TheAuthors}

\begin{document}

\maketitle

\begin{abstract}
  As CPU clock speeds have stagnated and high performance computers continue to have ever higher core counts, increased parallelism is needed to take advantage of these new architectures. Traditional serial time-marching schemes can be a significant bottleneck, as many types of simulations require large numbers of time-steps which must be computed sequentially. Parallel in Time schemes, such as the Multigrid Reduction in Time (MGRIT) method, remedy this by parallelizing across time-steps, and have shown promising results for parabolic problems. However, chaotic problems have proved more difficult, since chaotic initial value problems (IVPs) are inherently ill-conditioned. MGRIT relies on a hierarchy of successively coarser time-grids to iteratively correct the solution on the finest time-grid, but due to the nature of chaotic systems, small inaccuracies on the coarser levels can be greatly magnified and lead to poor coarse-grid corrections. Here we introduce a modified MGRIT algorithm based on an existing quadratically converging nonlinear extension to the multigrid Full Approximation Scheme (FAS), as well as a novel time-coarsening scheme. Together, these approaches better capture long-term chaotic behavior on coarse-grids and greatly improve convergence of MGRIT for chaotic IVPs. Further, we introduce a novel low memory variant of the algorithm for solving chaotic PDEs with MGRIT which not only solves the IVP, but also provides estimates for the unstable Lyapunov vectors of the system. We provide supporting numerical results for the Lorenz system and demonstrate parallel speedup for the chaotic Kuramoto-Sivashinsky partial differential equation over a significantly longer time-domain than in previous works.
\end{abstract}



\section{Introduction}

Although Parallel in Time (PinT) methods date back over 50 years~\cite{Ni1964}, interest in these methods has only recently increased due to the stagnation of CPU clock speeds beginning in the 2000s. See~\cite{gander_review,ben_jacob_review} for an introduction and overview of the field. For many problems, spatial parallelism can become saturated, while the time dimension remains largely unparallelized and solving for many time steps sequentially is still common practice. Although PinT schemes, if perfected, have great potential for speedup when combined with existing spatially parallel techniques, they have yet to see widespread adoption, because the time dimension presents difficulties not seen in the spatial dimension. The most important difficulty is \emph{causality,} obeying a one-way dependence where the solution at a later time depends solely on the solution at previous times, while the spatial dimension has coupling in multiple directions. To date, PinT has already been demonstrated to provide substantial speedups for parabolic problems~\cite{ben_jacob_review, parabolic_mgrit}, such as the heat equation, largely due to the fact that the causality of the system is relaxed over time. Parabolic problems have \emph{weak} dependence on initial conditions, since they tend toward a steady state which is largely uncorrelated to the initial data. Hyperbolic problems, such as the wave equation, remain difficult to parallelize in time, since they have \emph{strong} dependence on initial conditions, although some speedup has been demonstrated for such problems in special cases~\cite{ben_jacob_review, oliver_advection}. However, to our knowledge, no speedup has yet been achieved for chaotic problems, which exhibit \emph{sensitive} dependence on the initial conditions, with the result that the initial value problem (IVP) is ill-conditioned for chaotic systems. Despite these difficulties, chaotic systems are a very important class of problems with a wide range of applications across science and engineering, such as turbulent fluid flows, combustion, and climate modeling \cite{wang2013drag, kuramoto1978diffusion, lorenz1984irregularity}.

As a general trend, PinT methods are more attractive for problems involving many time-steps. This is because for smaller problems, the overhead introduced by many PinT methods prevent them from competing with sequential time-stepping. However, the ill-conditioning associated with chaotic problems cause PinT methods to converge slower as the time-domain length increases \cite{LSS_PinT}. Thus for many chaotic problems, there is no crossover point after which existing PinT methods are faster than sequential time-stepping. The difficulties of applying PinT methods to chaotic problems have been identified in several works. For example, \cite{PFASST} introduces the PFASST algorithm, a PinT method with robust convergence for advection diffusion equations. However, the authors observed a degradation in convergence for the chaotic Kuramoto-Sivashinsky (KS) equation, even over a time-scale 180 times smaller than the largest considered in this paper for the KS equation. Noteably, \cite{LSS_PinT} introduces a parallel in time algorithm which converges over arbitrary time-scales for chaotic problems, however, this is achieved by solving the linearized Least Squares Shadowing (LSS) problem, which relaxes the initial condition and thus yields the solution to a \emph{nearby} IVP, not the original one posed. Further, LSS, while well-conditioned, forms an expensive optimization problem, and requires a good initial guess over the entire time-domain.
Because of these limitations, the LSS method was not demonstrated to provide a speedup over sequential time-stepping for the Lorenz system nor the KS equation. 
Thus the principal challenge which we take on in this paper is to extend the time-domain length over which PinT methods can converge for chaotic problems, allowing for the first time a speedup over sequential time-stepping.

While there are many other promising PinT approaches, direct and iterative \cite{gander_review,ben_jacob_review}, we consider the Multigrid Reduction in Time (MGRIT)~\cite{MGRIT14} framework, and propose modifications that enable parallel speedup for chaotic IVPs. In particular, we investigate a modification to the the Full Approximation Scheme (FAS) coarse-grid equation used by  MGRIT, as well as propose a modified coarsening scheme for time-stepping propagators, which together greatly improve MGRIT convergence for the considered chaotic systems. The proposed modified FAS algorithm is equivalent to applying the Multilevel Nonlinear Method (MNM) algorithm~\cite{irad_MNM} (an iterative method which solves nonlinear \emph{elliptic} equations with quadratic convergence) to the \emph{time-dimension}. Thus the proposed algorithm yields a PinT method which we demonstrate numerically to have quadratic convergence. We further introduce a novel low-memory and matrix-free variant of the algorithm based on approximations of the Lyapunov vectors. This allows for fast convergence to the IVP solution while simultaneously generating estimates for the subset of backward Lyapunov vectors, which form a basis for the unstable manifold of the system. Thus the proposed algorithm provides the solution to the IVP, as well as estimates of the Lyapunov vectors, with a parallel-in-time speedup over just solving the IVP with sequential time-stepping.

This paper is organized as follows: We first introduce the standard MGRIT algorithm in Section \ref{sec:intro_MGRIT}, and study its poor performance on the chaotic Lorenz system in Section \ref{sec:intro_chaos}. Sections \ref{sub:theta} and \ref{sub:delta} motivate and introduce a new time-coarsening scheme (called the $\theta$ method) and the modified FAS coarse-grid correction (called $\Delta$ correction).  Section \ref{sub:lowrank} then introduces a novel memory efficient version of $\Delta$ correction based on approximations of the backward Lyapunov vectors.  Section \ref{sec:numerics} first demonstrates significantly improved MGRIT convergence for the chaotic Lorenz system and then a parallel speedup for a larger chaotic problem, the KS equation, through strong and weak scaling studies.  In both cases, we are careful to choose long time-domains designed to allow for significant chaotic behavior.  
Finally, Section \ref{sec:conclusions} concludes and proposes future work.

\subsection{MGRIT}%
\label{sec:intro_MGRIT}
MGRIT is an iterative multigrid method for solving discrete IVPs of the form
\begin{equation}
  \begin{cases}
    \vu_0 = \vg_0, \\
    \vu_{i+1} = \Phi(\vu_i) + \vg_{i+1} & i = 0, 1, 2, \dots, n-1,
  \end{cases}
  \label{eqn:ivp}
\end{equation}
where $\Phi$ is a nonlinear time-stepping operator. Systems of this form typically arise from a time-discretization of an ordinary differential equation (ODE) of the form $\vu'(t) = \vf(\vu(t))$, in which case $\Phi$ is some time-stepping scheme such as Euler's method and $\vg_i$ corresponds to constant forcing terms. The system is defined over a discrete time-grid with $N_t + 1$ points, $\fg = {\{t_i\}}_{i=0}^{N_t}$, and time-step size $h = t_{i+1} - t_i$. We will assume without loss of generality, 
that $h$ is constant. Let $\vu = \begin{bmatrix} \vu_0,&\vu_1,&\dots,&\vu_{N_t}\ \end{bmatrix}^T$ denote the state vector and let $\vg = \begin{bmatrix} \vg_0, & \vg_1, & \dots, & \vg_{N_t} \end{bmatrix}^T$ be a constant forcing term which also encodes the initial condition. Then equation~\eqref{eqn:ivp} may be written in the form of a block non-linear matrix equation,
\begin{equation}
  A (\vu) = \vg \text{, where } A(\vu) = \begin{bmatrix}
    I      &       &        &        &   \\
    - \Phi & I     &        &        &   \\
           & -\Phi & I      &        &   \\
           &       & \ddots & \ddots &   \\
           &       &        & - \Phi & I
  \end{bmatrix}
  \begin{bmatrix}
     \vu_0 \\ \vu_1 \\ \vu_2 \\ \vdots \\ \vu_{N_t}
  \end{bmatrix}
  \label{eqn:systemA}.
\end{equation}
 Typically, this system would be solved using forward substitution, which corresponds with traditional sequential time-marching. MGRIT instead applies FAS multigrid iterations to the system~\eqref{eqn:systemA}, allowing it to be solved in parallel. To this end, \eqref{eqn:systemA} is approximated on a hierarchy of coarser time grids, e.g. $\Omega^{2h}, \Omega^{4h}, \Omega^{8h}, \dots$, which provide error corrections to the finer grids, while the finer grids provide further corrections via local block Jacobi relaxation.

\begin{figure}[t]
  \centering
  \begin{tikzpicture}[xscale = 1.2]
    \draw [thick] (-4.5,0) -- (5.0,0);

    \draw[line width=2pt] (-4.5,-0.2) -- (-4.5,0.2);
    \draw (-4.0,-0.2) -- (-4.0,0.2);
    \draw (-3.5,-0.2) -- (-3.5,0.2);
    \draw (-3.0,-0.2) -- (-3.0,0.2);
    \draw (-2.5,-0.2) -- (-2.5,0.2);
    \draw[line width=2pt] (-2.0,-0.2) -- (-2.0,0.2);
    \draw (-1.5,-0.2) -- (-1.5,0.2);
    \draw (-1.0,-0.2) -- (-1.0,0.2);
    \draw (-0.5,-0.2) -- (-0.5,0.2);
    \draw (0,-0.2) -- (0,0.2);
    \draw[line width=2pt] (0.5,-0.2) -- (0.5,0.2);
    \draw (1.0,-0.2) -- (1.0,0.2);
    \draw (1.5,-0.2) -- (1.5,0.2);
    \draw (2.0,-0.2) -- (2.0,0.2);
    \draw (2.5,-0.2) -- (2.5,0.2);
    \draw[line width=2pt] (3.0,-0.2) -- (3.0,0.2);
    \draw (3.5,-0.2) -- (3.5,0.2);
    \draw (4.0,-0.2) -- (4.0,0.2);
    \draw (4.5,-0.2) -- (4.5,0.2);
    \draw[line width=2pt] (5.0, -0.2) -- (5.0, 0.2);

    \node [below] at (-4.5, -0.2) {$t_0$};
    \node [above] at (-4.5, 0.2) {$T_0$};
    \node [below] at (-4.0, -0.2) {$t_1$};
    \node [below] at (-3.5, -0.2) {$t_2$};
    \node [below] at (-3.0, -0.2) {$t_3$};
    \node [below] at (-2.5, -0.2) {$\cdots$};
    \node [below] at (-2.0, -0.2) {$t_m$};
    \node [above] at (-2.0, 0.2) {$T_1$};
    \node [below] at (5.0, -0.2) {$t_{N_t}$};
    \node [above] at (5.0, 0.2) {$T_{N_T}$};

    \draw [thick,decorate,decoration={brace,amplitude=6pt,raise=0pt}] (0.5, 0.2) -- (3.0, 0.2);
    \node [above] at (1.8, 0.5) {$m h$};
    \draw [thick,decorate,decoration={brace,amplitude=3pt,raise=0pt, mirror}] (1.0, -0.2) -- (1.5, -0.2);
    \node [below] at (1.3, -0.3) {$h$};
  \end{tikzpicture}
  \caption{Uniformly spaced fine-grid points and coarse-grid points with coarsening factor $m$. The $T_i$ are the C-points and form the coarse-grid, while the small hashmarks $t_i$  are F-points. Together, the F- and C-points form the fine-grid.\label{fig:grids}}
\end{figure}
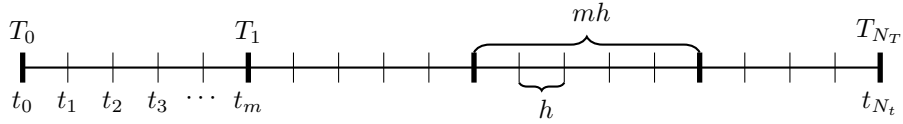

The multigrid method requires a coarsening scheme in time, inter-grid transfer operators, and a coarse-grid equation. Here, we define those for a two-level MGRIT method with fine grid, $\fg$, and coarse-grid, $\cg$, for coarsening factor $m$.  
To coarsen in time, label every $m$th time-point in $\fg$ a C-point and all other points an F-point, then $\cg$ is the set of size $N_T$ containing the C-points in $\fg$ (see Figure~\ref{fig:grids}). A C-point, along with the following $m-1$ F-points to the right, is called a coarse interval. For grid transfer operations, MGRIT uses \emph{injection}. For restriction, injection  maps the values of $\vu$ at the C-points in $\fg$ to the corresponding points in $\cg$, and for interpolation, it maps the points in $\cg$ to the corresponding C-points in $\fg$. The action of restriction by injection is given by the block matrix
\begin{equation}\label{eq:injection}
    R = \begin{bmatrix}
      I &   &        &   &   &        &   &        &   & \\
        & 0 & \cdots & 0 & I &        &   &        &   & \\
        &   &        &   &   & \ddots &   &        &   & \\
        &   &        &   &   &        & 0 & \cdots & 0 & I
    \end{bmatrix},
\end{equation}
yielding a coarse-grid vector $\vu_c = R \vu$ on $\cg$. 
Interpolation by injection is hence given by the action of $R^T$.
Following interpolation from $\cg$ to $\fg$, i.e. $\vu \gets R^T \vu_c$, MGRIT always relaxes the solution on $\fg$ using F-relaxation~\eqref{eqn:frelax}, which evolves the state at each C-point to the following F-points in each coarse interval using $\Phi$. 
F-relaxation may be viewed as part of the interpolation process, in that it updates each F-point based on the new C-point information.

The two-level MGRIT method as described thus far is equivalent to the popular Parareal algorithm \cite{GaVa2007}. However, apart from being a multi-level method, one of MGRIT's differences in comparison to Parareal is support for FCF-relaxation. Here F-relaxation is followed by C-relaxation~\eqref{eqn:crelax}, where the solution is propagated from the last F-point in each coarse interval to the following C-point, which is then followed by another F-relaxation~\cite{MGRIT14}. FCF relaxation has been shown in many cases to greatly improve convergence, especially in the multilevel setting~\cite{MGRIT14,Do2016}. Importantly, since the coarse-grid intervals are disjoint, F- and FCF- relaxation can be done in parallel, with minimal communication between processors. 
\begin{align}
    \vu_{km + i} &= \Phi(\vu_{km + i - 1}) + \vg_{km+i} &\text{for $i = 1, \dots, m-1$ and each coarse-interval $k$}\label{eqn:frelax} \\
    \vu_{km} &= \Phi(\vu_{km - 1}) + \vg_{km} &\text{for each coarse-interval $k$}\label{eqn:crelax}
\end{align}

Coarsening in time induces a new system of equations posed on $\cg$, where the ideal space-time operator $A_*$ is given by 
\begin{align}
    A_* = \begin{bmatrix}
    I        &         &        &          &   \\
    - \Phi^m & I       &        &          &   \\
             & -\Phi^m & I      &          &   \\
             &         & \ddots & \ddots   &   \\
             &         &        & - \Phi^m & I
  \end{bmatrix} 
\end{align}
 Here $\Phi^m(\vu_i)$ is understood to denote a fine-grid propagation of the solution across one coarse interval, from one C-point to the next, including the  forcing term $\vg$, e.g. $\Phi^2(\vv_{i-2}) = \Phi(\Phi(\vv_{i-2}) + \vg_{i-1})$, $\Phi^3(\vv_{i-3}) = \Phi(\Phi(\Phi(\vv_{i-3}) + \vg_{i-2}) + \vg_{i-1})$, etc. While solving this ideal coarse-grid equation yields the exact solution for each C-point, this would require as much work as solving the original fine-grid problem. Instead, MGRIT introduces the coarse-grid time-stepping operator $\Phi_c$ to cheaply approximate the action of $\Phi^m$ by rediscretizing the fine-grid equations on the coarse-grid, and then solves the coarse-grid equation with the approximate operator $A_c$: 
\begin{equation}%
  \label{eqn:splitting}
  A_* \approx A_c =
  \begin{bmatrix}
    I        &         &        &          &   \\
    - \Phi_c & I       &        &          &   \\
             & -\Phi_c & I      &          &   \\
             &         & \ddots & \ddots   &   \\
             &         &        & - \Phi_c & I
  \end{bmatrix}.
\end{equation}

In the case that $\Phi$ is derived from a continuous-time problem, then $\Phi_c$ is usually derived from some (possibly nonstandard) rediscretization of the continuous problem over the new coarse time-grid $\cg$. Deriving coarse operators for general $\Phi$ is an open problem~\cite{ben_jacob_review,hyperbolic_coarse_grid} and motivates this paper.

In the case that $R$ is injection and interpolation is followed by an exact F-relax, the approximation of $A_*$ by $A_c$ may be interpreted as a non-linear splitting method. Let $\vec{\tau}(\vu_c) = A_c(\vu_c) - A_*(\vu_c)$, then $A_*(\vu_c) = A_c(\vu_c) - \vtau(\vu_c) = \vg_c$, and one immediately gets the well-known $\tau$-correction form of FAS multigrid~\cite{achi_77}: 
\begin{equation}%
  \label{eqn:FAS_coarse_grid}
  A_c(\vv_c^{k+1}) = \vg_c + \vec{\tau}(\vv_c^{k}),
\end{equation}
where $\vv_c^k$ denotes an approximate coarse solution on $\cg$ after $k$ multigrid iterations, and for MGRIT, $\vec{\tau}_i = \Phi^m(\vv^k_{i-1}) - \Phi_c(\vv^k_{i-1})$. This splits the operator $A_*$ into a part that is cheap to invert, $A_c$, and a part which can be computed efficiently in parallel, $\vtau$. One iteration of the two-level MGRIT scheme involves computing $\vtau(\vv^k)$ at the C-points on $\fg$, injecting $\vg$ and $\vtau$ to $\cg$, solving~\eqref{eqn:FAS_coarse_grid} sequentially for $\vv^{k+1}$, and then interpolating to $\fg$ and applying F-relaxation.
The vector $\vtau$ takes the form of a forcing term on the coarse-grid, and it steers the solution toward the fine-grid solution, as well as ensuring that the exact fine-grid solution is a fixed point of the iteration.
To see this, we plug the exact solution $\vu$ satisfying $A_*(R \vu) = \vg_c$ into~\eqref{eqn:FAS_coarse_grid}:
\begin{align*}
  A_c(\vv_c^{k+1})     & = \vg_c + \vtau(R \vu)      \\
                       & = \vg_c + A_c(R \vu) - A_*(R \vu) \\  
                       & = A_c(\vu_c) \\
  \implies \vv_c^{k+1} & = R \vu.
\end{align*}
Thus any $\vu$ satisfying the ideal coarse-grid equation is a fixed point of the MGRIT iteration, which is not true in general without the use of the $\vtau$ correction on the coarse-grid.
The two-level MGRIT algorithm is detailed in Algorithm~\ref{alg:MGRIT}.
A multi-level MGRIT algorithm then results from recursive application of the two-level scheme to solve the coarse-grid equation~\eqref{eqn:FAS_coarse_grid} with another MGRIT cycle.
\newcommand{\MGRIT}[1]{$\mathrm{MGRIT}_{#1}$} 
Recursion of this process gives the V-cycle \MGRIT{m_l} algorithm, with $m_l$ time-grid levels.

\begin{algorithm}%
  \caption{MGRIT two grid cycle: \MGRIT{2}$(\vv, \vg, m)$}
  \begin{algorithmic}
    \FOR{each C-point, $i = 1, 2, 3, \dots, N_T$}
    \STATE $\vec{\tau}_i \gets \Phi^m(\vv_{c,i-1}) - \Phi_c(\vv_{c,i-1})$
    \ENDFOR
    \STATE{$\vg_c \gets R \vg$ and solve:}
    \FOR{$i = 1, 2, 3, \dots, N_T$}
    \STATE{$\vv_{c,i} \gets \Phi_c(\vv_{c,i-1}) + \vg_{c,i} + \vec{\tau}_i$}
    \ENDFOR{}
    \STATE{refine, then relax with $\Phi$ (F-relax or FCF-relax)}
  \end{algorithmic}
  \label{alg:MGRIT}
\end{algorithm}

\subsection{Motivation: Chaotic problems and MGRIT}%
\label{sec:intro_chaos}

To study MGRIT for chaotic systems, we will use the Lorenz system as a model problem. The Lorenz system is a three-dimensional system of ODEs which is widely studied as an archetypal example of a chaotic system, and is given by
\begin{equation}
  \begin{cases}
    x' & = \sigma (y-x)    \\
    y' & = x(\rho - z) - y \\
    z' & = xy - \beta z
  \end{cases}.
  \label{eqn:lorenz}
\end{equation}
For the classical values of parameters $\sigma = 28$, $\rho = 10$, and $\beta = 8/3$, the Lorenz system is chaotic, with greatest \emph{Lyapunov exponent} (LE) of $\lambda_1 \approx 0.9$~\cite{strogatz}. This can be understood to mean that two trajectories differing only infinitesimally in initial conditions will, almost surely, diverge exponentially from each other in time with average rate $\lambda_1$. Although trajectories diverge from each other, the Lorenz system is \emph{Lyapunov stable}, meaning that trajectories are ultimately confined to a bounded trapping region in space, which contains a \emph{strange attractor}, a fractal manifold which is the limit set of the Lorenz system. Generally, a system with $n_s$ spatial dimensions has $n_s$ LEs, which are characteristic of the qualitative behavior of the system, and every chaotic system has a greatest LE which is greater than zero.
The corresponding \emph{Lyapunov vectors} (LVs) are characteristic directions, $\psi^k(t)$, $k=1, 2, \dots, n_s$, along which infinitesimal perturbations will grow exponentially with average rate $\lambda_k$. For example, the Lorenz system is three-dimensional and has three LEs: $\lambda_1 \approx 0.9$ which corresponds with perturbations lying tangent to the surface of the strange attractor, $\lambda_2 = 0$ which corresponds with perturbations tangent to the flow (resulting in a difference only in phase), and $\lambda_3 \approx -14$, which corresponds with perturbations away from the strange attractor. The set of LVs having negative LEs is referred to as the \emph{stable manifold}, the set having vanishing LEs is the \emph{neutral manifold}, and the set having positive LEs is the \emph{unstable manifold}.

PinT simulations of chaotic systems such as Lorenz are difficult because of two main problems. The first is that the sensitivity of chaotic systems makes them ill-conditioned. The other is that coarsening in time can cause serious global qualitative changes in the behavior of the discretized system because of changes in the Lyapunov spectrum between the fine- and coarse-grid. While chaotic systems are best known for their sensitivity to perturbations in the initial conditions, they are equally sensitive to changes to the system parameters. Coarsening in time can be considered a parametric perturbation, and thus the solution on the naive coarse-grid will diverge exponentially fast from that of the fine-grid. Therefore, the challenge is to form a coarse-grid equation that is both locally precise and that also captures the global qualitative behavior of the system, even for very coarse time-grids.

\subsubsection{Conditioning and halting criterion}
Due to the existence of the unstable manifold, any numerical error committed in solving a chaotic IVP with any method will grow exponentially in time and eventually result in error which is of the same order of magnitude as the solution itself. We can characterize this with the condition number of the IVP~\eqref{eqn:ivp}. First, we define \emph{Lyapunov time}, $T_\lambda = \frac{\ln(10)}{\lambda_1}$, to be the time it takes for a perturbation to a trajectory to grow by a factor of 10~\cite{strogatz}. Hence, $T_\lambda$ is the time it takes for our numerical simulation of the system to lose one digit of accuracy due to the ill conditioning. 
The  Lyapunov time can be used to compute an estimate for the condition number, $\kappa = \mathcal{O}(10^{T_f/T_\lambda})$, where $T_f = N_t h$ is the length of the time-domain.  Lyapunov time serves as a normalized time-scale over which all chaotic systems have greatest Lyapunov exponent equal to unity.  For example, one Lyapunov time for a weather simulation might be on the order of a couple of days, while one Lyapunov time for a simulation of planetary motion might be on the order of millions of years.  Thus, our consideration of time domains of length 4 and 8 in Lyapunov time represents an interval of significant length.

Traditionally, MGRIT halts based on a space-time residual and that is the strategy pursued here.  This strategy ignores the forward error, which measures the difference between the MGRIT and sequential time-stepping solutions, grows exponentially along the unstable manifold, and can be large.  Instead, for ill-conditioned chaotic IVPs, it is appropriate to consider \emph{backward} error. Here, we consider the Shadowing Lemma \cite{Pilyugin1999}, which provides backward error guarantees for discrete, numerical trajectories over a wide class of systems, such as the Lorenz system.  In such cases, when the forward error is large, as long as an approximate trajectory has a \emph{residual}, here given by $\vr_i = \vg_i + \Phi(\vu_{i-1}) - \vu_i$, which is small and uniformly bounded, there exists an exact solution to a nearby problem which is uniformly close to the approximate trajectory.  That is, the backward error is small.  Thus our residual convergence criteria depends on this heuristic, where as long as our PinT method has reduced the global residual below a small uniform tolerance, we assume the Shadowing Lemma applies and that the obtained solution has backward error which is small and uniformly bounded. See \cite{what_good_are, Pilyugin1999} for more details about solving chaotic IVPs with numerical methods.

\subsubsection{Poor naive MGRIT performance for chaotic problems}
\begin{figure}[h]
  \centering
  \includegraphics[width=\textwidth]{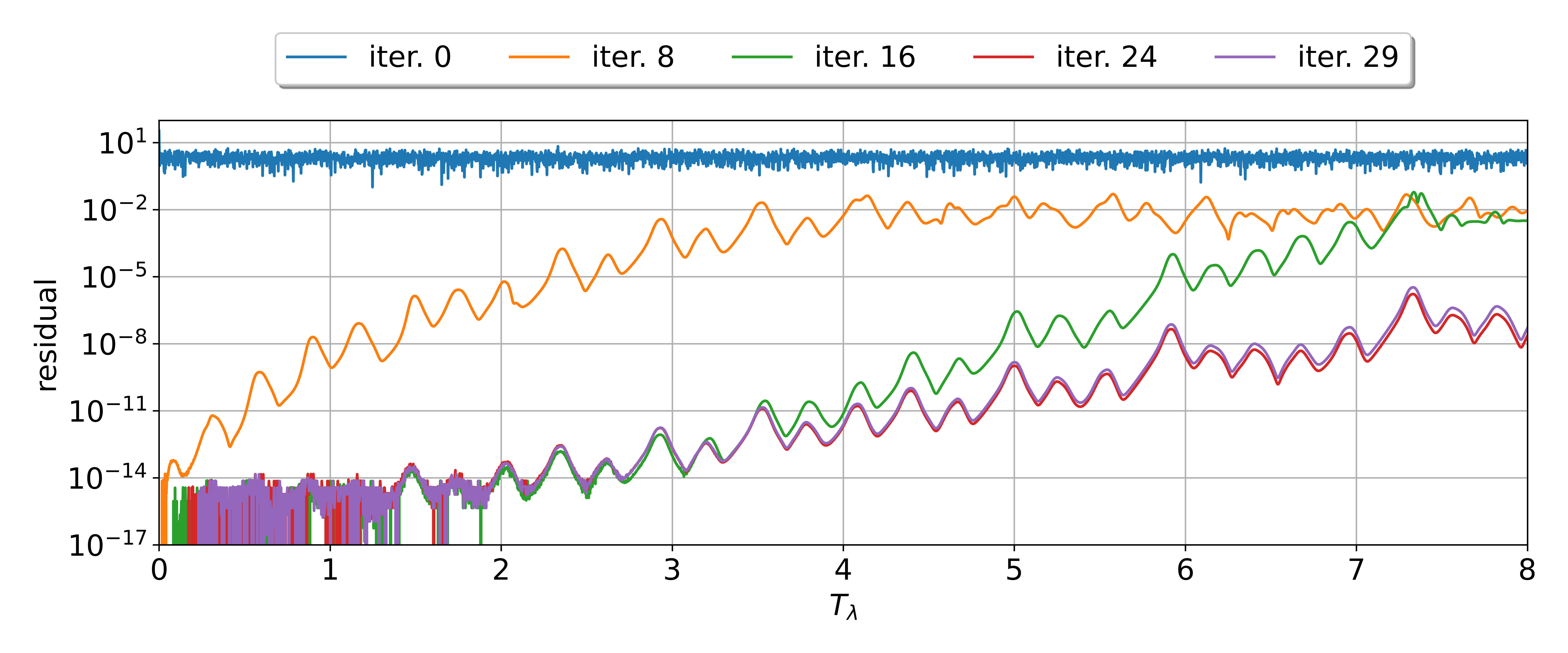}
  \caption{Residual over the time-domain (in units of Lyapunov time) for 30 two-level MGRIT iterations on the Lorenz equation with $T_f = 8T_\lambda$ and coarsening factor $m=2$.  Convergence stalls after 24 iterations, due to the exponentially growing residual which is not damped by further iterations.\label{fig:MGRIT_stall}}
\end{figure}

For the Lorenz system, we choose a final time of $T_f = 8T_\lambda$, giving the condition number $\kappa = \mathcal{O}(10^8)$, and exemplify MGRIT performance in Figure~\ref{fig:MGRIT_stall}. This figure demonstrates that the residual grows exponentially in time, with average rate $\lambda_1$, and thus the solution after 30 iterations does not satisfy our backward error argument, since its residual is large at the later times and MGRIT iterations have stagnated. It is clear that for chaotic problems, MGRIT alone does not suffice to converge to an acceptable solution within a reasonable number of iterations. 

\section{Main contributions}%
\label{sec:main}

\subsection{Improved time-coarsening scheme: \texorpdfstring{$\theta$}{Theta} method}%
\begin{figure}[h]
  \centering
  \includegraphics[width=0.8\textwidth]{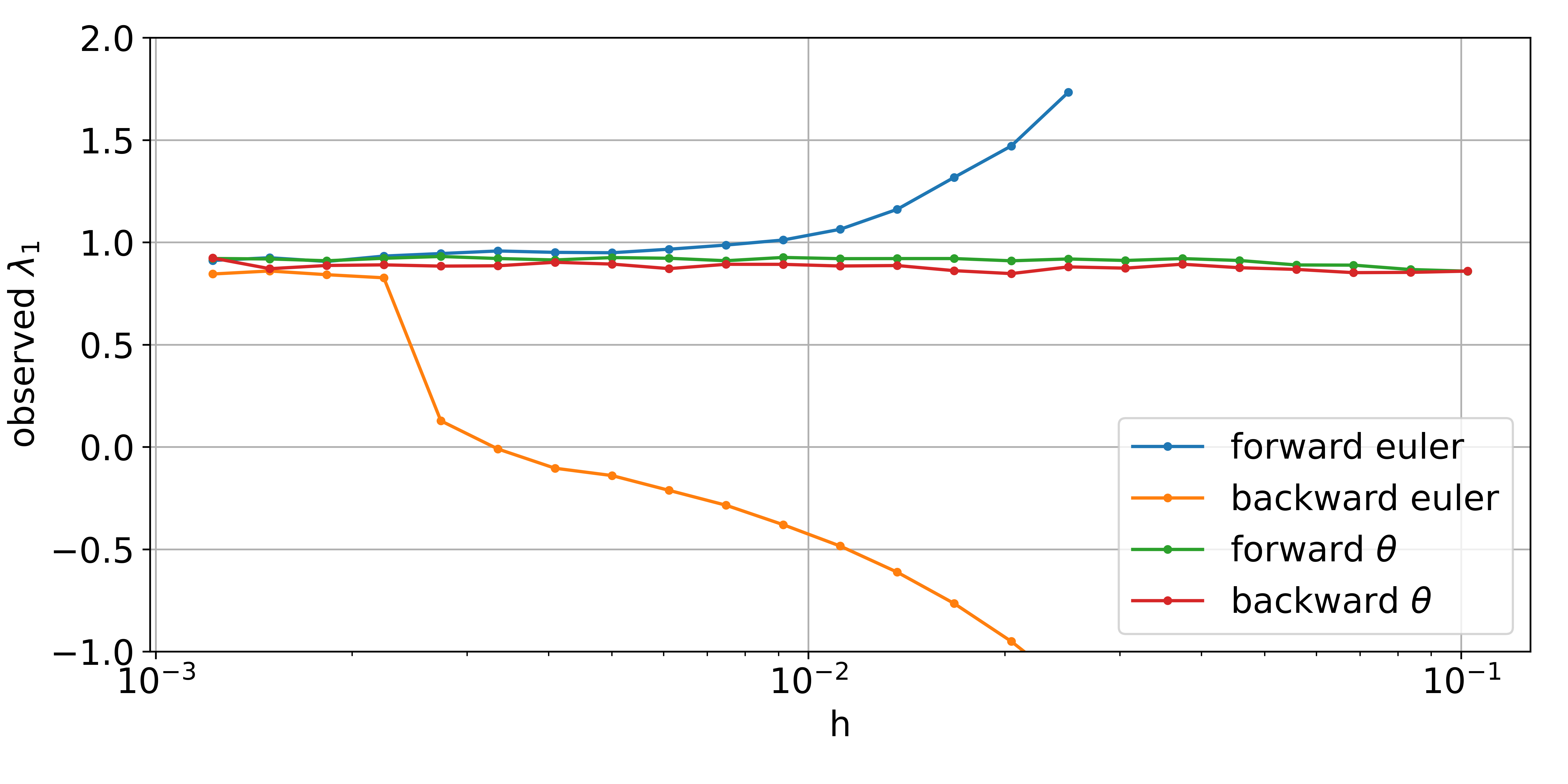}
  \caption{Plot of the observed greatest Lyapunov exponent $\lambda_1$ for different time-step sizes $h$, comparing forward Euler, backward Euler, as well as the $\theta$ method proposed here. For forward and backward Euler, coarsening in time changes the qualitative behavior of the system, while the $\theta$ method, using the ``forward" and ``backward" values of $\theta$, preserves this sensitivity on coarse-grids.}
  \label{fig:lyap_dependence}
\end{figure}

\label{sub:theta}
One difficulty in solving chaotic systems with PinT is that coarsening in time can cause dramatic qualitative changes to the global behavior of the system.  Typically for MGRIT, when an explicit method is used on the fine grid, an implicit method is used on very coarse-grids for stability. However, switching from an explicit method to an implicit one causes serious changes to the behavior of the discretized chaotic system. For example, when using backward Euler to solve the Lorenz equations, the measured greatest Lyapunov exponent decreases with increasing time-step size, $h$, meaning that for large $h$, a chaotic system can become artificially stabilized. Conversely, using forward Euler, the Lyapunov exponent increases with increasing $h$, and the system appears \emph{more} chaotic on coarse-grids~\cite{what_good_are}. 
Thus, a time-stepping scheme that preserves the qualitative behavior of the system on coarse-grids, is one which lies somewhere between the binary of forward and backward Euler.  Figure~\ref{fig:lyap_dependence} demonstrates this dependence for different time-stepping schemes applied to the Lorenz system, including for the $\theta$ method proposed here.

As shown by the previous example, a suitable fine-grid time-stepping scheme does not necessarily make a suitable coarse-grid scheme. One option would be to use a higher order method on the coarse-grid, thus ensuring better accuracy on coarse-grids. However, higher order methods are generally more expensive, and in our experiments do not improve MGRIT convergence.  Although higher order methods approximate the \emph{continuous} problem to a higher order, this is not useful, since the coarse-grid operator instead needs to approximate the discrete fine-grid problem. In fact, a higher order $\Phi_c$ will only agree with $\Phi^m$ to the same order of accuracy as $\Phi$ agrees with the continuous problem, so in a naive implementation, the accuracy of the approximation $\Phi_c \approx \Phi^m$ is limited by the lowest order of accuracy among the pair $\Phi_c, \Phi$, meaning that it is not useful in this context to use a more expensive, higher order method on the coarse-grid.  Thus, we propose the $\theta$ method, which does not refer to a specific time-discretization, but is rather a framework for choosing a coarse-grid propagator which better approximates the discrete fine-grid solution than simple rediscretization.  Here, we rely on the fact that the fine-grid time-step size, $h$, is small enough that the time-stepping operators $\Phi_c$ and $\Phi^m$ admit convergent Taylor series in $h$ and can therefore be classified by their formal orders of accuracy.  The more general case, where either $h$ or $m$ is large enough that this does not hold, will be explored as future work.

The $\theta$ method relies on three principles: First, the coarse-grid time-stepper, $\Phi_c$, should have at least the same order of accuracy as the fine-grid, $\Phi$, so that we have a suitably consistent coarse-grid discretization. Second, rather than choosing a single $\Phi_c$, we propose to consider a $k$-parameter family of coarse-grid time-steppers, $\Phi_c(\vu_i, \vec{\theta})$, whose parameters $\theta_j$, $j = 1, 2, \dots, k$ can be tuned depending on the physics of the problem and the time-step size, $h$.  Third, when possible, these parameters are tuned such that $\Phi_c$ approximates $\Phi^m$ to a higher order than $\Phi$ approximates the continuous equation. The choice of the parameterization itself depends on the physics of the IVP and will typically blend explicit and implicit time-stepping schemes.

More precisely, consider a fine-grid time stepping operator of order $p$ with $m$-step stability function $\phi_m(z)$, representing $m$ consecutive fine-grid steps.  We propose to consider a $k$-parameter family of time-stepping operators of at least order $p$, with stability function $\phi_\theta(mz)$, representing a single large coarse-grid step. If the system of equations
\begin{subequations}
\begin{align}
  \phi_m(0)       & = \phi_\theta(0)       \\
  \phi_m'(0)      & = \phi_\theta'(0)      \\
                  & \vdots                 \nonumber \\
  \phi_m^{(p+1)}(0) & = \phi_\theta^{(p+1)}(0) \\
                  & \vdots                 \nonumber \\
  \phi_m^{(p+k)}(0) & = \phi_\theta^{(p+k)}(0)
\end{align}
\end{subequations}
has a solution for $\{\theta_j\}_{j=1}^k$, then these parameters yield a method whose stability function $\phi_\theta$ approximates the $m$-step stability function $\phi_m$ to order $p + k$. 

We exemplify this idea by deriving a coarse-grid propagator $\Phi_c$ to better approximate $m$ steps of forward Euler, as in the example above.  Let the discretized ODE be given by $u_t = f(u)$, then this problem can be solved to first order by any member of the family of implicit first order single step methods given by 
\begin{equation}
  \vu_{i+1} = \vu_i + h [\theta f(\vu_i) + (1 - \theta)f(\vu_{i+1})],
  \label{eqn:theta}
\end{equation}
where the parameter $\theta \in [0, 1]$ controls the explicit/implicit balance. 

Since the $\theta$ method is first order for \emph{any} value of $\theta$, we can use this extra degree of freedom to better approximate the fine grid \textit{discretization} by choosing $\theta$ in such a way that $m$ fine-grid steps of forward (or backward) Euler are approximated to \emph{second} order in $mh$. For example, apply forward Euler and the $\theta$ method to the scalar Dalquist problem, $u' = \lambda u$ for complex constant $\lambda$, and let $z = h\lambda$, where $h$ is the time-step size. Then, the stability functions for $m$ steps of forward Euler with time-step $h$ and for the $\theta$ method with time-step $mh$ are given by
\begin{align}
   \phi_m(z)         = \phi(z)^m = (1 + z)^m\;\; \mbox{ and }\;\;   
  \phi_\theta(z)  = \frac{1 + \theta mz}{1 - (1 - \theta) mz},
\end{align}
respectively. In order for $\phi_\theta$ to approximate $\phi_m$ to second order in $z$, the following three equations must be satisfied:
\begin{align}
  \phi_m(0)    = \phi_\theta(0),    \quad
  \phi_m'(0)   = \phi_\theta'(0),   \;\; \mbox{ and }\;\; 
  \phi_m''(0)  = \phi_\theta''(0).
\end{align}
While the first and second equality are already satisfied since forward Euler and the $\theta$ method are both first order in $z$, the third equality yields
\begin{subequations}
\begin{align}
  m(m - 1)     & = 2(1 - \theta_m)m^2 \\
  1 - \theta_m & = \frac{m - 1}{2m}   \\
   \theta_m     & = \frac{m + 1}{2m}. \label{eqn:theta_ass}
\end{align}
\end{subequations}
Since these values of $\theta$ give a method which approximates the fine-grid to 2nd order in the scalar case, the natural question is whether this method also better represents the system on coarse-grids for nonlinear, multivariable problems. This is indeed the case, as depicted in Figure \ref{fig:lyap_dependence}, where the $\theta$ method is compared to forward and backward Euler applied to the Lorenz system. 
The $\theta$ method accurately preserves the greatest Lyapunov exponent even on coarse time-grids. Thus for this example, we can say that the $\theta$ method addresses, at least in part, the fundamental difficulties with PinT and chaotic problems, preserving the global dynamics on coarse time grids and improving local accuracy.

In section~\ref{sec:numerics_ks}, we discuss a second order $\theta$ method which we have constructed for the solution of the KS equation, a stiff, nonlinear, chaotic PDE, and demonstrate the improvement to MGRIT performance that it provides relative to rediscretization.

\subsection{Improved FAS coarse-grid equation: \texorpdfstring{$\Delta$}{Delta} correction}\label{sub:delta}

As we have seen, naive MGRIT fails to converge for the Lorenz system over long time-domains due to a residual which grows exponentially in time.  For systems where the underlying dynamics are chaotic, MGRIT is very sensitive to errors, no matter how small. This is because while the $\vec{\tau}$ correction makes the fine-grid solution a fixed point on the coarse-grid, there may still be a significant mismatch between the Lyapunov spectra on the coarse and fine grids.  Such a mismatch is why even small errors on the coarse-grid are exponentially magnified along incorrect Lyapunov coordinate vectors, causing the observed exponentially increasing residual.  This mismatch can be measured (and eventually corrected for) by considering the tangent linear propagator, $F_i$, along the trajectories of the fine and coarse operators.  We will see that $F_i$ on the coarse-grid needs to be extremely accurate relative to the fine-grid.

Since we know that perturbations only grow along the unstable manifold for trajectories of the Lorenz system, we should expect that the components of the error along the unstable manifold are the slowest to converge, while the other components of the error converge more quickly.  We further do not want the error in the unstable manifold, which may be large, to affect the residual.  However in order to better understand how chaos effects the MGRIT cycle, we need a more concrete definition for the LVs, which form a basis for the Lyapunov manifolds of the Lorenz system, and we need to understand how perturbations along the Lyapunov manifolds are propagated in finite time. Thus, it is helpful to establish a distinction between \emph{forward} and \emph{backward} LVs, $\psi^{+, k}(t)$ and $\psi^{-, k}(t)$, respectively.

Let $\vu(t_i)$ be a discrete trajectory of a nonlinear dynamical system having time propagator $\Phi$, which passes through some point $\vu(t_0) = \vu_0$. Then, the linearization of the time-stepper $\Phi$ at each point $\vu_i$ yields the matrix $F_i = D_u \Phi(\vu_i)$, where $D_u$ is the differential operator with respect to the spatial variables. The matrix $F_i$ is called the tangent linear propagator \cite{lyap_vecs}, as it describes the propagation of infinitesimal perturbations along $\vu$, i.e. $\Phi(\vu_i + \ve) - \vu_{i+1} \approx F_i \ve$. Thus, the propagation of a perturbation from time point $t_i$ to point $t_j$ is given by $W(t_i, t_j) = F_j F_{j-1} \dots F_{i+1} F_i$. In the limit as $j \to \infty$, the time-average of the singular values of $W(t_i, t_j)$ are equal to $\exp(\lambda_i)$, where the $\lambda_i$ are the LEs of the system. The LEs are independent of the times $t_i$, $t_j$, and are also the same for almost all $\vu_0 \in \mathbf{R}^{n_s}$, and are thus considered constants of the system. The forward LVs at time $t_i$, $\psi^{+, k}(t_i)$, are given by the time-average of the right singular vectors of $W(t_i, t_j)$ in the limit as $j\to\infty$, while the backward LVs at time $t_j$, $\psi^{-, k}(t_j)$, are given by the time-average of the left singular vectors in the limit as $i \to -\infty$. Note that while the LEs are not time-dependent, the backward and forward LVs are.  For brevity, we will use the phrase ``unstable manifold" to refer to the subspace which is spanned by the backward LVs having positive LE, although this is not in a strict sense the typical definition of the unstable manifold, which requires the definition of covariant LVs that are beyond the scope of this paper.

\subsubsection{Understanding Lyapunov vectors and MGRIT}

Assume that the LEs are distinct and ordered from greatest to least, i.e. $j > k \implies \lambda_j < \lambda_k$, and let $\Psi^-(t_i)$ be the orthonormal matrix whose columns are the ordered backward LVs at time $t_i$. The forward propagation of $\Psi^-$ in finite time is given by \cite{lyap_vecs}:
\begin{equation}\label{eq:lyap_QR_relation}
    F_i \Psi^-(t_i) = \Psi^-(t_{i+1}) R_{i+1},
\end{equation}
where $R_{i+1}$ is upper triangular.  Thus, equation \eqref{eq:lyap_QR_relation} indicates that, forward in time, the \emph{backward} LVs are mapped to vectors which are orthogonal to the set of LVs with smaller LE.  For the Lorenz system, a perturbation along the unstable manifold will remain orthogonal to the neutral and stable manifolds for all time, while a perturbation along the stable manifold will almost surely have a nonzero projection on the neutral and unstable manifolds in finite time.

\begin{figure}[h]
    \centering
\begin{subfigure}{.5\textwidth}
  \centering
  \includegraphics[width=\linewidth]{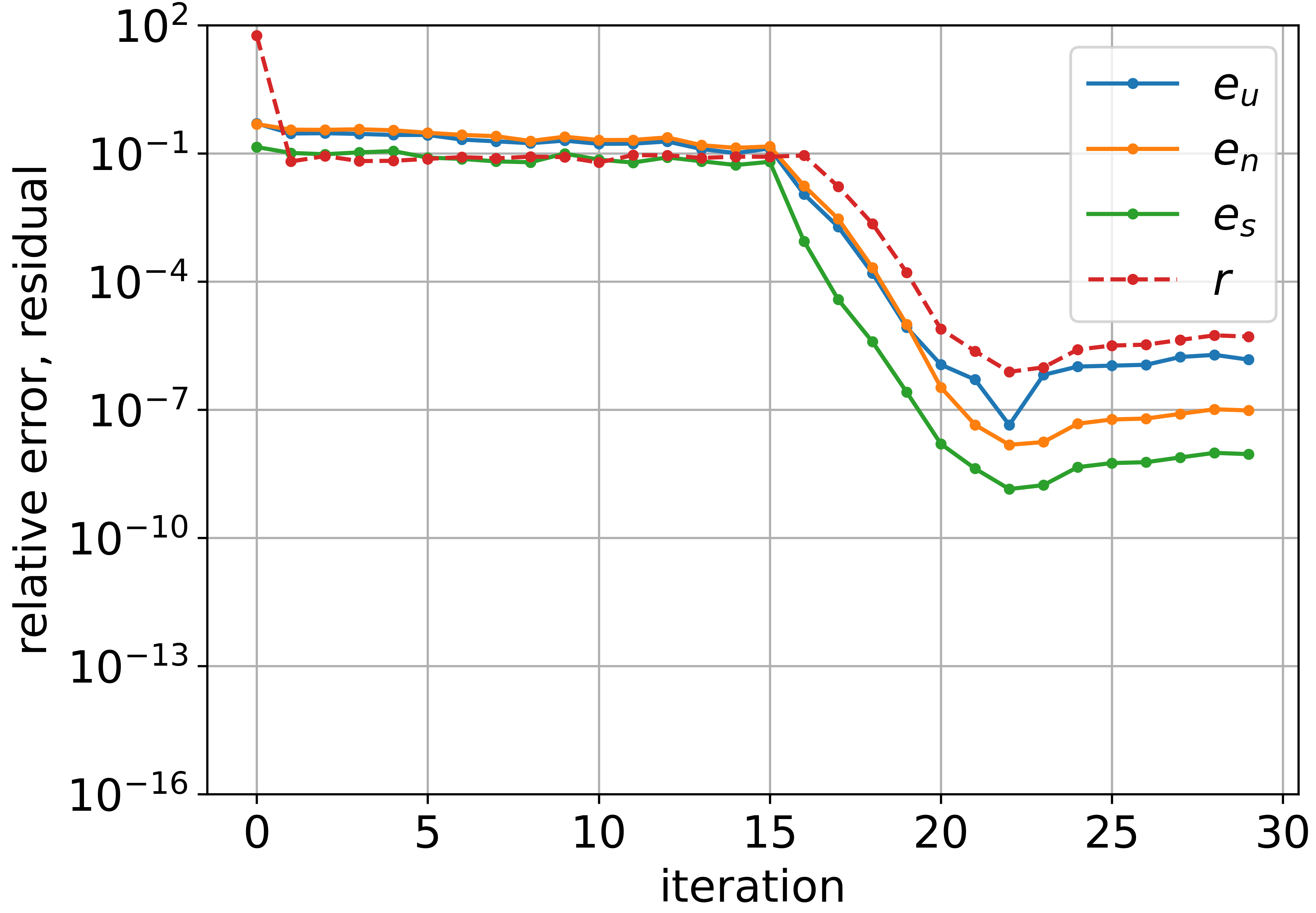}
  \caption{Classical MGRIT.\label{fig:stall_decomp_MGRIT}}
  \label{fig:sub1}
\end{subfigure}%
\begin{subfigure}{.5\textwidth}
  \centering
  \includegraphics[width=\linewidth]{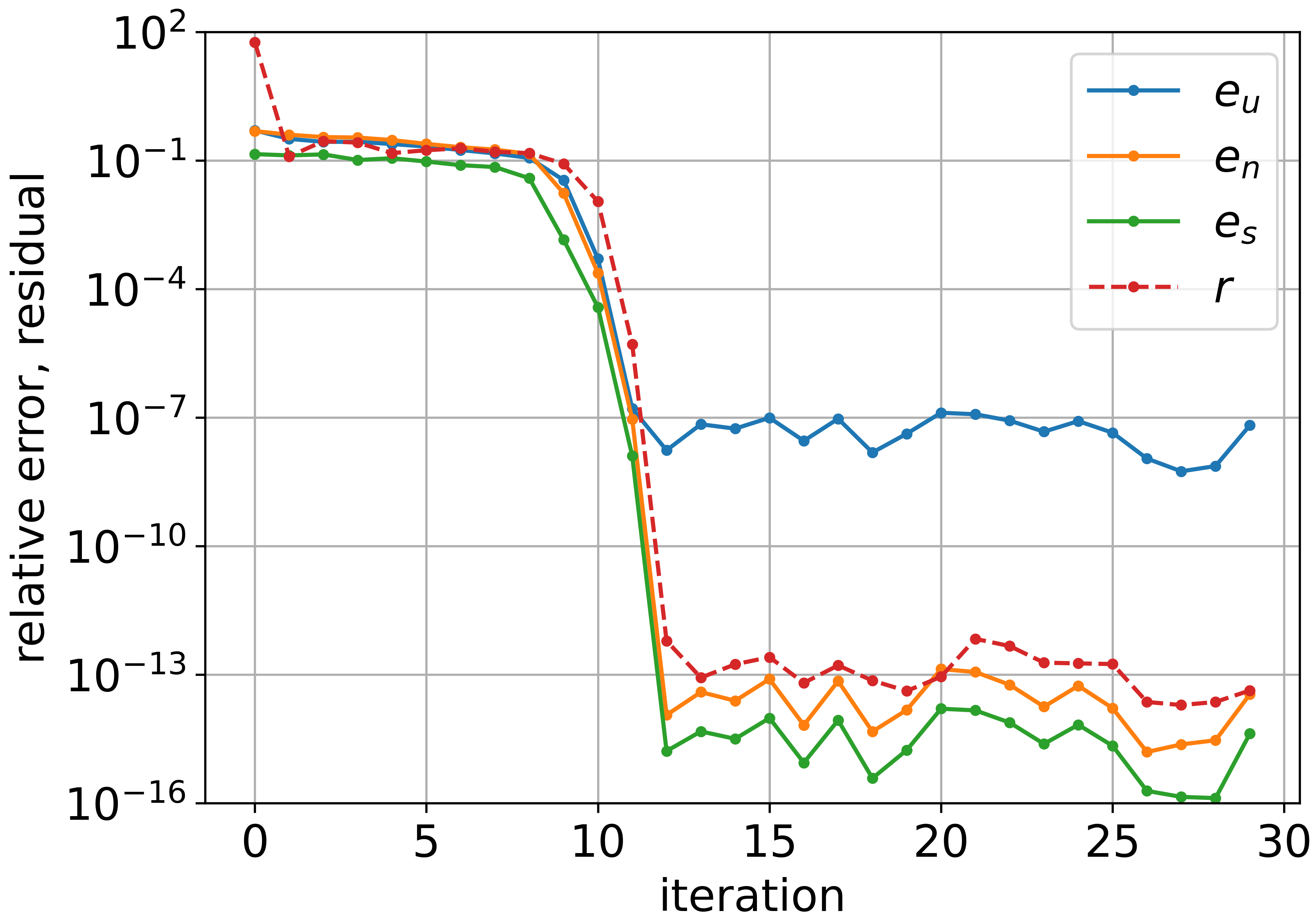}
  \caption{MGRIT with $\Delta$ correction.\label{fig:stall_decomp_Delta}}
  \label{fig:sub2}
\end{subfigure}
\caption{The 2-norm of the residual ($r$) and the components of relative error along the unstable, neutral, and stable manifolds, ($e_u$, $e_n$, $e_s$, respectively) for 30 two-level MGRIT iterations on the Lorenz equation with $T_f = 8T_\lambda$ and coarsening factor $m=2$.  (\ref{fig:stall_decomp_MGRIT}) Convergence stalls after 22 iterations, and there is significant error even in the stable manifold, which MGRIT should damp. (\ref{fig:stall_decomp_Delta}) The Delta correction prevents $e_u$ from affecting the residual, and the algorithm is able to converge despite significant error in the unstable manifold. \label{fig:stall_decomp}}
\end{figure}

Returning to MGRIT, assuming there is a large component of error along the \emph{backward} unstable manifold, we would expect that error to remain orthogonal to the other manifolds, and if there is a small component of error along the stable backward manifold, it should only contribute a small amount to the error in the less stable manifolds. Thus, if MGRIT is capable of damping errors along the stable manifold but not the unstable, we should be able to observe this by computing the components of error using the matrices $\Psi^-(t_i)^T$ as a time-dependent change of basis. However, in practice, since the coarse-grid uses an approximate time-stepping operator, $\Phi_c$, the coarse-grid equation will not have the same Lyapunov spectrum. This mismatch between the LVs on the fine and coarse-grid cause some of the error along the unstable backward manifold to ``leak" into the neutral and stable manifolds, thus stalling convergence. Figure \ref{fig:stall_decomp_MGRIT} demonstrates this phenomenon, that after some initial convergence, the error in the unstable manifold contributes significantly to the error in the other manifolds as well, preventing the residual from decreasing further.  

\subsubsection{Deriving \texorpdfstring{$\Delta$}{Delta} correction}
The $\Delta$ correction remedies the mismatch in the Lyapunov spectrum between the fine and coarse-grids by using \emph{linearizations} of the fine operator, computed about the current solution guess, to correct the coarse operator.

Let $\Phi$, $\Phi_c$, and $\Phi^m$ be defined as above, and $\vv = \vu - \ve$ be an approximate solution on the fine-grid. Now define
\begin{equation}
  \Delta_i := \left(D_u \Phi^m - D_u \Phi_c\right)(\vv_i),
  \label{eqn:delta}
\end{equation}
where $\Delta_i$ is a matrix valued function of $\vv_i$, which encodes the difference between the linearizations of the ideal and coarse operators. 
The matrix $\Delta_i$ will naturally have the same number of dimensions as the number of spatial dimensions of the system.
Contrast this with the $\vec{\tau}$ correction, which encodes the difference between the \emph{values} of these two operators applied to $\vv_i$.
We then use the computed $\Delta_i$ to form a correction to the time-stepper on each coarse interval:
\begin{equation}
  \Phi_{\Delta_i}(\vv_{c,i}) := \Phi_c(\vv_i) + \Delta_i \vv_{c,i},
  \label{eqn:phi_coarse}
\end{equation}
which ensures that as $\vv$ approaches $\vu$, i.e.\ near MGRIT convergence, $D_u \Phi_\Delta$ approaches $D_u \Phi^m$.  Because $D_u \Phi_\Delta$ is the linear tangent propagator on the coarse grid, which determines the Lyapunov spectrum, the $\Delta$ correction is able to correct the mismatch in the Lyapunov spectrum, even when $\Phi_c$ is a poor approximation to $\Phi^m$.

Together with the $\vec{\tau}$ correction, which is computed at the same time, this gives the modified MGRIT algorithm \ref{alg:delta}, where the new additions are colored in red.

\begin{algorithm}
  \caption{MGRIT two-level cycle with $\Delta$ correction: $\Delta \mathrm{MGRIT}_2(\vv, \vg, m)$\label{alg:delta}}
  \begin{algorithmic}
    \FOR{each C-point, $i = 1, 2, 3, \dots, N_T$}
    \STATE \textcolor{red}{$\Delta_i \gets D_u \Phi^m(\vv_{c,i-1}) - D_u \Phi_c(\vv_{c,i-1})$}
    \STATE $\vec{\tau}_i \gets \Phi^m(\vv_{c,i-1}) - \Phi_{\Delta_i}(\vv_{c,i-1})$
    \ENDFOR
    \STATE{$\vg_c \gets R \vg$ and solve:}
    \FOR{$i = 1, 2, 3, \dots, N_T$}
     \STATE{$\vv_{c,i} \gets $\textcolor{red}{$\;\Phi_{\Delta_i}(\vv_{c,i-1})$} $+\;\vec{\tau}_i + \vg_{c,i}$}
    \ENDFOR
    \STATE{interpolate, then relax on $\vv$ with $\Phi$ (F- or FCF-relax)}
  \end{algorithmic}
\end{algorithm}

Note that the first loop does not update the values of $\vec{v}$ at each time point, and may thus be done in parallel, while the loop on the coarse-grid must be solved sequentially. Remember that the multilevel method replaces the sequential solve on the coarse-grid with a recursive call to the algorithm. In Figure \ref{fig:stall_decomp_Delta}, we see the effect this $\Delta$ correction has on the convergence of MGRIT for the Lorenz system. Without the $\Delta$ correction, Figure \ref{fig:stall_decomp_MGRIT} shows that there are significant components of error along the neutral and stable manifolds, causing a stall in residual convergence, however, the $\Delta$ correction in Figure \ref{fig:stall_decomp_Delta} completely eliminates this problem, allowing the iteration to converge in residual. Thus we see that the $\Delta$ correction addresses the main difficulty with applying MGRIT to chaotic systems by correcting the mismatch in the Lyapunov spectrum between the fine and coarse-grids.

\subsubsection{Modified FAS coarse-grid equation and quadratic convergence}

Since the $\Delta$ correction updates the coarse operator, this represents a modification to the FAS coarse-grid equation. Recall that in Section~\ref{sec:intro_MGRIT} we saw that the splitting $\vec{\tau}(\vv_c) = A_c(\vv_c) - A_*(\vv_c)$ resulted in the FAS coarse-grid equation~\eqref{eqn:FAS_coarse_grid}. Now, define the global linear correction on the coarse-grid as $\Delta(\vv_c) = D_u A_*(\vv_c) - D_u A_c(\vv_c)$, which in turn implies a modified coarse-grid equation of
\begin{equation}
  \label{eqn:modified_FAS}
  [A_c + \Delta(\vv_c^k)](\vv_c^{k+1}) = \vg_c + \vec{\tau}(\vv_c^{k}) + [\Delta(\vv_c^k)]\vv_c^k.
\end{equation}
Equation \eqref{eqn:modified_FAS} is equivalent to the Multilevel Nonlinear Method (MNM) from \cite{irad_MNM} applied to the time dimension. The MNM was shown numerically to converge quadratically for nonlinear elliptic PDEs, thus we expect our $\Delta$ corrected MGRIT algorithm to exhibit quadratic convergence as well.

For the simplest choice of coarse-grid propagator where $\Phi_c \equiv 0$, we have that $A_c = I$, the identity operator, and~\eqref{eqn:modified_FAS} becomes 
\begin{align*}
  [I + D_u A_*(\vv^k) - I](\vv^{k+1}) & = \vg_c + \vec{\tau}(\vv^k) + [D_u A_*(\vv^k) - I](\vv^k)    \\
  \implies [D_u A_*(\vv^k)]\vv^{k+1}           & = \vg_c + \vv^k - A_*(\vv^k) + [D_u A_*(\vv^k)]\vv^k - \vv^k \\
 \implies [D_u A_*(\vv^k)]\vv^{k+1}           & = [D_u A_*(\vv^k)]\vv^k + \vg_c - A_*(\vv^k)                 \\
  \implies \vv^{k+1}                  & = \vv^k - [D_u A_*(\vv^k)]^{-1}(A_*(\vv^k) - \vg_c),
\end{align*}
which is equivalent to an iteration of Newton's method applied to the residual equation $A_*(\vv^k) - \vf_c$. We should then expect MGRIT with this $\Delta$ correction to converge at least as well as Newton's method, i.e. local quadratic convergence, as long as $A_c$ approximates $A_*$ better than the identity $I$. In Section \ref{sec:numerics}, we provide numerical evidence for the quadratic convergence of MGRIT with $\Delta$ correction.

\subsection{Low-rank \texorpdfstring{$\Delta$}{Delta} correction for PDEs}
\label{sub:lowrank}
For the Lorenz system with dimension $n_s=3$, the computation and storage requirements for the $\Delta$ correction are small relative to the improvement in convergence they provide. However, the computation and storage of $n_s \times n_s$ matrices will be prohibitive for the case of most PDEs, which often have $n_s$ in the millions or larger.
However, it is often the case that only a small, finite number of the LEs are positive, meaning that the unstable manifold of the discretized system has dimension much smaller than $n_s$.  Moreover, there is evidence that the number of positive LEs can be bounded in many cases. For example in \cite{turbulence_spectrum}, the authors measure the Lyapunov spectrum for a 2D simulation of a chaotic flow around an airfoil at high Reynolds number while refining in space. They observed that although the Lyapunov spectrum of the system was mesh dependent, beyond a certain resolution threshold, the total number of unstable modes were no more than 5, even as the number of degrees of freedom were increased to $\mathcal{O}(10^6)$.

Since the negative Lyapunov exponents behave similar to a parabolic problem, we should expect MGRIT to converge well along these directions \cite{Do2016}, even without $\Delta$ correction. In contrast, we have already seen that the unstable modes must be represented very accurately on the coarse-grid, or else exponentially growing error will be mapped incorrectly to the stable manifold. This motivates a low-rank approximation to the $\Delta$ correction which targets only the  low-dimensional unstable manifold. Given an approximate solution $\{\vu_i\}_{i=0}^{n_t}$, let $\Psi_i$ be a rectangular, orthonormal matrix with columns equal to the first $k$ backward Lyapunov vectors $\psi_i^j$ for $j=0,1,\dots k-1$, then
\begin{equation}\label{eqn:lr_Delta}
  \hat{\Delta}_i = \Delta_i \Psi_i \Psi_i^T
\end{equation}
is a rank $k$ approximation to $\Delta_i$, which is exact for the first $k$ backward LVs.
Next, note that since $\Delta_i$ is a linearization of the function $\Phi^m - \Phi_c$, the columns of $\Delta_i \Psi_i = (D_u \Phi^m - D_u \Phi_c) \cdot \Psi_i = D_{\Psi_i} \Phi^m - D_{\Psi_i} \Phi_c$ are just directional derivatives along the first $k$ LVs. Thus, we only need to compute these $k$ directional derivatives and $\Psi_i$ in order to form $\hat{\Delta}_i$, never forming the full matrix $\Delta_i$. This requires storage of the factors $\Delta_i \Psi_i$ and $\Psi_i$, which are $n_s \times k$ matrices, meaning that as long as $k < n_s/2$, the low rank approximation requires less storage and less computational work than the full $\Delta$ correction.

However, implementation of this approach requires overcoming the cost of computing the first $k$ Lyapunov vectors, which are considered expensive to obtain and are classically computed sequential in time. The expense is largely due to the time propagation of LV matrices and linearizations of the time-stepping operator, as well as the QR factorizations needed for orthogonalization. Here, we propose an efficient way to estimate the LVs simultaneously while solving the original equation parallel-in-time. Given a trajectory $\{\vu_i\}_{i=0}^{\infty}$ and an arbitrary initial orthogonal $n_s \times k$ matrix $\Psi_0$, $k \leq n_s$, the QR iteration
\begin{equation}%
  \label{eqn:QR_iteration}
   \Psi_{i+1} R_i = [D_u \Phi(\vu_i)]\Psi_i \quad i = 0, 1, \dots
\end{equation}
will result in convergence of the columns of $\Psi_i$ to the first $k$ backward LVs of the trajectory as $i \to \infty$. For this reason, the backward LVs are sometimes called Gram-Schmidt vectors. The long-time average of the diagonals of the upper triangular matrix $R_i$ will converge to $e^{h\lambda_j}$, where $\lambda_j$ are the LEs. The QR decomposition orthonormalizes the columns of $\Psi_i$, preventing numerical instability, however in practice, $\Psi_i$ may be normalized only every $m$ steps while still converging to the LVs.  For details on computing LVs, see the work~\cite{lyap_vecs}. 

Given a finite trajectory, this QR iteration will yield an estimate to the true LVs which is more accurate toward the end of the time-domain. However, the sequential computation of these LV estimates is often more expensive than solving the state equation sequentially, thus we propose to compute these estimates parallel in time using MGRIT. Since the LV computation takes the form of an IVP, we may apply MGRIT directly to~\eqref{eqn:QR_iteration}, using $D_u \Phi_c$ as the coarse-grid time-stepping operator. In this way, our MGRIT cycle will be simultanesouly solving for the state vector $\vu$ and the LVs $\Psi_i$. It is important to note that while the IVP for chaotic systems is very sensitive to initial conditions, this is \emph{not} the case for the LVs. The first column of $\Psi_i$ will almost surely converge to the first backward Lyapunov vector $\psi^{-,1}(t_i)$, independent of the initial matrix $\Psi_0$. This closely resembles the convergence of the linear power iteration to the first eigenvector of a matrix. Therefore, although the state equation is chaotic, and thus ill-conditioned, the orthonormalization process ensures that the IVP for the backward LVs is well-conditioned, and thus we observe fast MGRIT convergence for the LVs without modification of the algorithm, and perhaps counter-intuitively, the most unstable modes will also be the most accurate, as they are the fastest to converge in the QR iteration algorithm.

Regarding computational cost, if we only normalize $\Psi_i$ with $QR$ at C-points and use \(\Psi_{i+1} = [D_u \phi(\vu)]\Psi_i\) at F-points, then we have the potential to save a great deal of computation without much loss in accuracy of the computed LVs.

In order for MGRIT to solve for the LVs $\Psi_i$, we also need to consider the MGRIT FAS coarse-grid for the problem of finding LVs. Let $F^{m} = D_u \Phi^m(\vu_i)$ and $F_c = D_u \Phi_c(\vu_i)$, then an appropriate $\tau$ correction term on the coarse-grid for the LVs is given by
\[
  \vtau_i = (F^{m} - F_c)\Psi_i.
\]
However, if we are using the low-rank $\Delta$ correction, then this becomes
\begin{align*}
  \vtau_i &= (F^{m} - (F_c + \Delta \Psi_i \Psi_i^T))\Psi_i \\
          &= (F^{m} \Psi_i - F_c \Psi_i - \Delta \Psi_i ) \\
          &= (\Delta - \Delta)\Psi \\
          &= 0.
\end{align*}
Thus as long as we are also using the low-rank $\Delta$ correction, the $\tau$ correction for the LVs can be ignored, as it is 0. From another perspective, the $\Delta$ correction is already a $\tau$ correction for the Lyapunov vectors.

We now have a low-rank $\Delta$ correction for the state equation and we have an efficient algorithm for finding the needed LV estimates given a trajectory. Combining the two, we get the low-rank $\Delta$ correction Algorithm~\ref{alg:lrDelta}, where changes relative to original MGRIT (Algorithm \ref{alg:MGRIT}) are again highlighted in red. In this algorithm, F-relaxation on the state vector is always followed by F-relaxation on the LVs, and likewise for C-relaxation. The LVs are also solved for sequentially on the coarsest grid alongside the state vector.  Modified Gram-Schmidt is used to compute the $QR$-factorization.

\begin{algorithm}
  \caption{MGRIT two grid cycle with LR$\Delta$ correction and LV estimates: $\Delta_{k} \mathrm{MGRIT}_2(\vv, \vg, m)$\label{alg:lrDelta}}
  \begin{algorithmic}
    \FOR{each C-point, $i = 1, 2, 3, \dots, N_T$}
    \STATE \textcolor{red}{$\hat{\Delta}_i \gets [D_u \Phi^m(\vv_{c,i-1}) - D_u \Phi_c(\vv_{c,i-1})] \Psi_{c,i-1}$}
    \STATE $\vec{\tau}_i \gets \Phi^m(\vv_{c,i-1}) - \Phi_{\Delta_i}(\vv_{c,i-1})$
    \ENDFOR
    \STATE{$\vg_c \gets R \vg$ and solve:}
    \FOR{$i = 1, 2, 3, \dots, N_T$}
    \STATE{$\vv_{c,i} \gets $\textcolor{red}{$\;\Phi_{\Delta_i}(\vv_{c,i-1})$} $+\;\vec{\tau}_i + \vg_{c,i}$}
     \STATE{\textcolor{red}{$\Psi_{c,i} \gets \mathrm{GramSchmidt}([D_u \Phi_{\Delta_i}(\vv_{c,i-1})]\Psi_{c,i-1})$}}
    \ENDFOR
     \STATE{interpolate, then relax on $\vv$ with $\Phi$, \textcolor{red}{and on $\Psi$ with $D_u \Phi$} (F- or FCF-relax)}
  \end{algorithmic}
\end{algorithm}

The parallel performance of MGRIT is highly dependent on the cost of the sequential solve on the coarse-grid, and the sequential propagation of the LVs on the coarse-grid adds significantly to this cost.  However, we have found that this sequential solve can either be skipped or approximated with a parallel relaxation step, without much change to the effectiveness of the algorithm in many cases. In fact, for many of the experiments presented in Section \ref{sec:numerics_ks}, no coarse-grid propagation of LVs is performed, since it was determined experimentally that FCF-relaxation on fine grids resolves the LV estimates sufficiently well for MGRIT convergence.  In this case, the only extra work is performed in parallel on the fine-grid (i.e., the GramSchmidt steps on the coarse-grid are skipped).  However for long time-domain sizes, the coarse-grid propagation of LVs is necessary in order to prevent the residual stalling, as discussed in the results section. 

\section{Numerical results}%
\label{sec:numerics}
We implemented the described algorithms by modifying the XBraid software package \cite{xbraid}, which is written in C and MPI. In all cases, the stopping criterion used is an absolute tolerance on the global 2-norm of the residual.

\subsection{The Lorenz system} %

In the following experiments, we solve the discret-ized Lorenz system~\eqref{eqn:lorenz} using forward Euler on the fine grid. A coarsening factor of $m=2$ is used across all of the studies. When the $\theta$ method \eqref{eqn:theta} is used on a coarse-grid, the values of $\theta$ are dependent on the grid level $l$ and are computed according to \eqref{eqn:theta_ass} with $m = 2^{l}$, where $l$ refers to the grid level, starting at $l=0$ for the fine-grid. The implicit equation \eqref{eqn:theta} is solved numerically using Newton's method. When the $\theta$ method is not used on the coarse-grid, forward Euler is used, with coarsened time-step size $m^{l}h$. Only F-relaxation is considered in this section. First, we study convergence rates for various two-level algorithms on a small problem. Then we perform a refinement study in time and a time-domain size scaling study. Finally we explore the effect of adding more coarse levels for different problem sizes.  

\subsubsection{Two-level results}

\begin{figure}[h]
  \centering
  \includegraphics[width=0.9\textwidth]{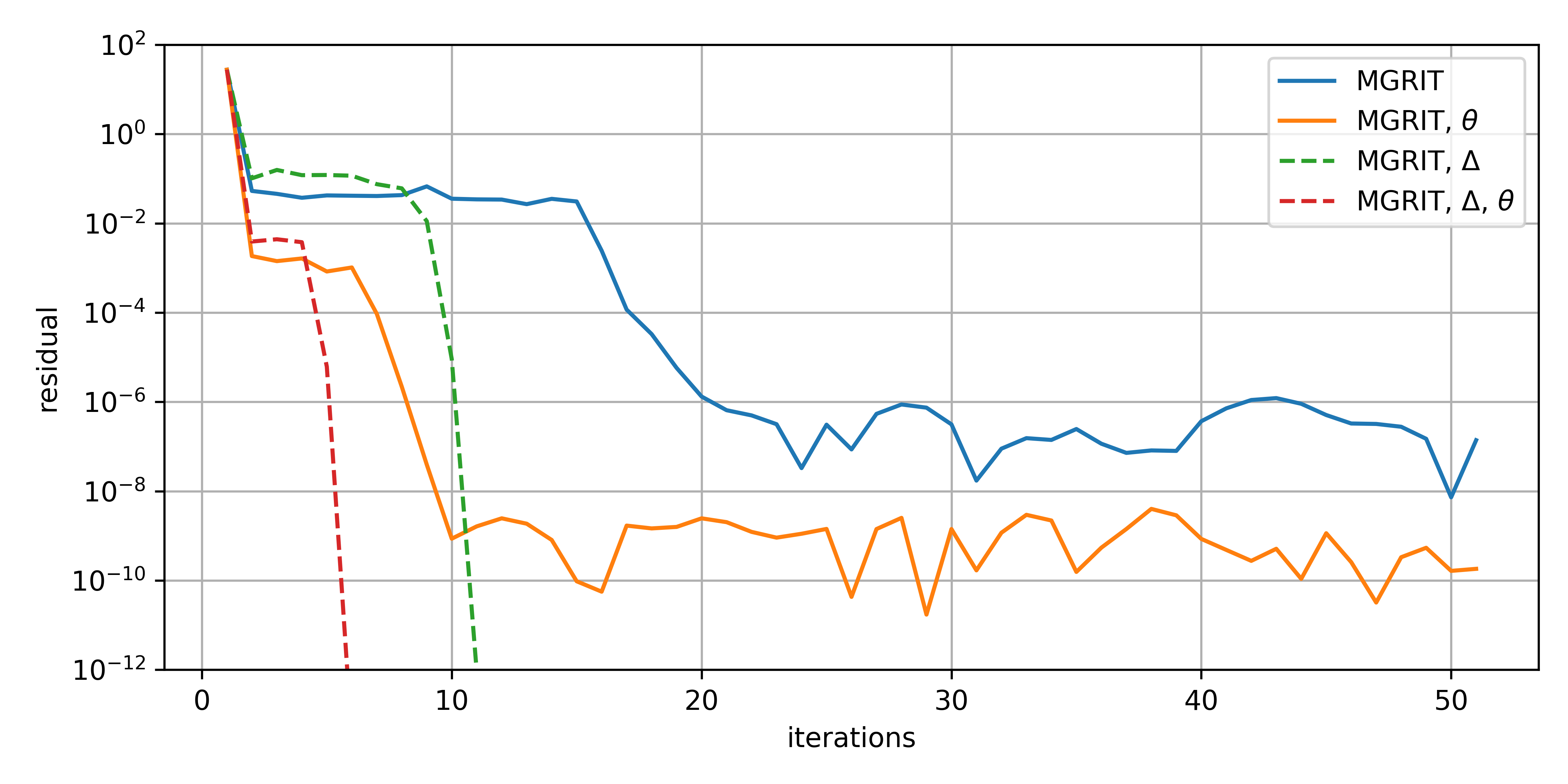}
  \caption{Residual history for each two-level algorithm applied to the Lorenz system with $T_f = 8T_\lambda$ and 8192 time-points. We see that standard \MGRIT{2} stalls on this time domain, as expected, but that we are able to converge quickly using the $\Delta$ correction. Using $\Delta$ correction and the $\theta$ method leads to the fastest convergence. \label{fig:convergence}}
\end{figure}

Figure~\ref{fig:convergence} plots the convergence history of the modified two-level MGRIT algorithms,
solving the Lorenz system with $T_f = 8T_\lambda$. This experiment demonstrates that \MGRIT{2}, even when using the $\theta$ method on coarse-grids, stalls for long enough time-domains, which is expected given the significant mismatch between the fine and coarse-grids regarding the Lyapunov spectrum. However, the $\Delta$ correction allows the method to converge quickly, with the combination of $\Delta$ correction and the $\theta$ method leading to the fastest convergence. Our intuition is that the increased accuracy of the $\theta$ method in approximating the fine-grid widens the basin of attraction for the quadratic convergence region of the $\Delta$ corrected method. 

In order to study the effect of varying time-step size $h$ on \MGRIT{2} performance, Table \ref{tab:refine_2} shows numbers of iterations required to reach a certain residual tolerance for the Lorenz system with fixed $T_f$ and increasing numbers of time-points $n_t$. In all cases fewer iterations are needed for smaller $h$.  However compared to naive \MGRIT{2}, both the $\Delta$ correction and the $\theta$ method require roughly half as many iterations to converge, and when used together, they require a quarter of the number of  iterations. Further, we see that the $\theta$ method can correct for instabilities on the coarse-grid which cause MGRIT to diverge (see first column).

\begin{table}
  \centering
  \caption{Iterations required for each two-level algorithm to converge to a residual tolerance of $10^{-10}$ for the Lorenz system with $T_f = 8T_\lambda$ and increasing number of time-points $n_t$. `*' indicates that the algorithm diverged.}
  \begin{tabular}{lccccc}
    \toprule
    \multicolumn{1}{c}{}                   & \multicolumn{5}{c}{$T_f$, $n_t$}                                         \\
    \cmidrule(rl){2-6}
    Algorithm
                                           & 4, 512                           & 4, 1024 & 4, 2048 & 4, 4096 & 4, 8192 \\
    \cmidrule(r){1-1} \cmidrule(rl){2-6}
    $\mathrm{MGRIT}_2$                     & *                                & 44      & 22      & 15      & 12      \\
    $\mathrm{MGRIT}_2$, $\theta$           & 19                               & 13      & 9       & 7       & 6       \\
    $\mathrm{MGRIT}_2$, $\Delta$           & *                                & 11      & 8       & 6       & 6       \\
    $\mathrm{MGRIT}_2$, $\Delta$, $\theta$ & 8                                & 6       & 5       & 4       & 4
    \\
    \bottomrule
  \end{tabular}
  \label{tab:refine_2}
\end{table}

Table \ref{tab:time_scaling_2} shows iteration counts for the two-level algorithm on the Lorenz system with increasing time-domain size $T_f$ and fixed time-step size $h$. For naive \MGRIT{2}, iteration counts increase linearly up until the critical time $T_f = 6T_\lambda$, after which naive \MGRIT{2} stalls. In contrast, the $\Delta$ correction and $\theta$ method greatly improve convergence for all time-domain sizes. Notably, the iteration counts for the $\Delta$ corrected algorithm are nearly flat (until the last column), even for long time-domain sizes, with the combined $\Delta$ method and $\theta$ method providing the fastest convergence. 

\begin{table}
  \centering
  \caption{Iterations required for each two-level algorithm to converge to a residual tolerance of $10^{-10}$ for the Lorenz system with varying $T_f$ (in Lyapunov time) and varying number of time-points $n_t$ such that $h$ is constant. `-' indicates that the algorithm did not converge within 100 iterations}
  \begin{tabular}{lcccccc}
    \toprule
    \multicolumn{1}{c}{}                   & \multicolumn{6}{c}{$T_f$, $n_t$}                                                         \\
    \cmidrule(rl){2-7}
    Algorithm
                                           & 2, 4096                          & 4, 8192 & 6, 12288 & 8, 16384 & 10, 20480 & 12, 24576 \\
    \cmidrule(r){1-1} \cmidrule(rl){2-7}
    $\mathrm{MGRIT}_2$                     & 10                               & 13      & 17       & 64       & -         & -         \\
    $\mathrm{MGRIT}_2$, $\theta$           & 4                                & 5       & 6        & 7        & 41        & -         \\
    $\mathrm{MGRIT}_2$, $\Delta$           & 5                                & 6       & 7        & 8        & 9         & 94        \\
    $\mathrm{MGRIT}_2$, $\Delta$, $\theta$ & 3                                & 4       & 4        & 5        & 5         & 48
    \\
    \bottomrule
  \end{tabular}
  \label{tab:time_scaling_2}
\end{table}

\subsubsection{Multilevel results}\label{sub:lorenz_ml}
While \MGRIT{2} is not typically used in practice, it is used as a stepping stone toward understanding the multilevel algorithm. Recall that \MGRIT{2} solves the coarse-grid equation \eqref{eqn:FAS_coarse_grid} using a sequential solve and that \MGRIT{3} replaces this direct sequential solve with a recursive application of \MGRIT{2} to inexactly solve the coarse-grid.  Thus, we expect that \MGRIT{3} will converge slower than \MGRIT{2}, and as we add more levels, this trend should continue. However, since the coarsest grid is solved sequentially, it is very important for parallel performance that the coarse-grid be as small as possible, since the proportion of the algorithm which is not parallelizable limits the scaling of the method.

Table \ref{tab:ml_scaling} depicts the effect of increasing the number of coarse-grids on MGRIT convergence for various algorithmic configurations.  For naive MGRIT ($\mathrm{MGRIT}_k$), the effects of coarsening beyond two levels quickly leads to divergence and an unusable method.  For the $\theta$ method ($\mathrm{MGRIT}_k$, $\theta$), convergence improves.  However for the 7 level solver (where the coarsest grid size is small enough to be practical), the method diverges or takes too many iterations (43) to be practical.  For $\Delta$ correction ($\mathrm{MGRIT}_k$, $\Delta$), convergence with 3 levels is impressively stable across all $T_f$, but for 5 and 7 levels divergence is observed.  This is primarily due to the fact that the time-stepping scheme becomes unstable on such coarse-grids. Finally, the combination of both approaches ($\mathrm{MGRIT}_k$, $\Delta$, $\theta$) combines the stable convergence for long $T_f$ with the coarse-grid stability of the $\theta$ method. The result is promising convergence at 7 levels for $T_f=2$ and $T_f=4$, since the coarsest grid here is small (64 and 128 time-points respectively) and the iteration counts are similar to previous cases demonstrating parallel speedup for linear parabolic problems \cite{MGRIT14}.  Unfortunately, the small spatial size of this problem (3) makes a parallel performance study difficult, as computations would always be dominated by network latency.  Thus, we next consider a larger problem.

\begin{table}
  \centering
  \caption{Iterations required for each algorithm, using varying numbers of grids, to converge to a residual tolerance of $10^{-10}$ for the Lorenz system with varying $T_f$ Lyapunov time and varying number of time-points $n_t$ such that $T_f/n_t$ is constant. `-' indicates that the algorithm did not converge within 100 iterations, `*' indicates that the algorithm diverged.}
  \begin{tabular}{lcccc}
    \toprule
    \multicolumn{1}{c}{}                   & \multicolumn{4}{c}{$T_f$, $n_t$}                                 \\
    \cmidrule(rl){2-5}
    Algorithm                              & 2, 4096                          & 4, 8192 & 6, 12288 & 8, 16384 \\
    \cmidrule(r){1-1} \cmidrule(rl){2-5}
    $\mathrm{MGRIT}_2$                     & 10                               & 13      & 17       & 64       \\
    $\mathrm{MGRIT}_3$                     & 13                               & 18      & -        & -        \\
    $\mathrm{MGRIT}_5$                     & 26                               & -       & -        & -        \\
    $\mathrm{MGRIT}_7$                     & *                                & *       & *        & *        \vspace{5pt} \\
    $\mathrm{MGRIT}_2$, $\theta$           & 4                                & 5       & 6        & 7        \\
    $\mathrm{MGRIT}_3$, $\theta$           & 6                                & 7       & 9        & 11       \\
    $\mathrm{MGRIT}_5$, $\theta$           & 10                               & 13      & 19       & 63       \\
    $\mathrm{MGRIT}_7$, $\theta$           & 43                               & -       & -        & -        \vspace{5pt}\\
    $\mathrm{MGRIT}_2$, $\Delta$           & 5                                & 6       & 7        & 8        \\
    $\mathrm{MGRIT}_3$, $\Delta$           & 6                                & 8       & 11       & 13       \\
    $\mathrm{MGRIT}_5$, $\Delta$           & *                                & *       & *        & *        \\
    $\mathrm{MGRIT}_7$, $\Delta$           & *                                & *       & *        & *        \vspace{5pt}\\
    $\mathrm{MGRIT}_2$, $\Delta$, $\theta$ & 3                                & 4       & 4        & 5        \\
    $\mathrm{MGRIT}_3$, $\Delta$, $\theta$ & 3                                & 4       & 5        & 5        \\
    $\mathrm{MGRIT}_5$, $\Delta$, $\theta$ & 5                                & 6       & 7        & 9        \\
    $\mathrm{MGRIT}_7$, $\Delta$, $\theta$ & 9                                & 15      & 20       & 23       \\
    \bottomrule
  \end{tabular}
  \label{tab:ml_scaling}
\end{table}

\subsection{The Kuramoto-Sivashinsky equation}%
\label{sec:numerics_ks}

\begin{figure}[h]
    \centering
    \includegraphics[width=0.8\textwidth]{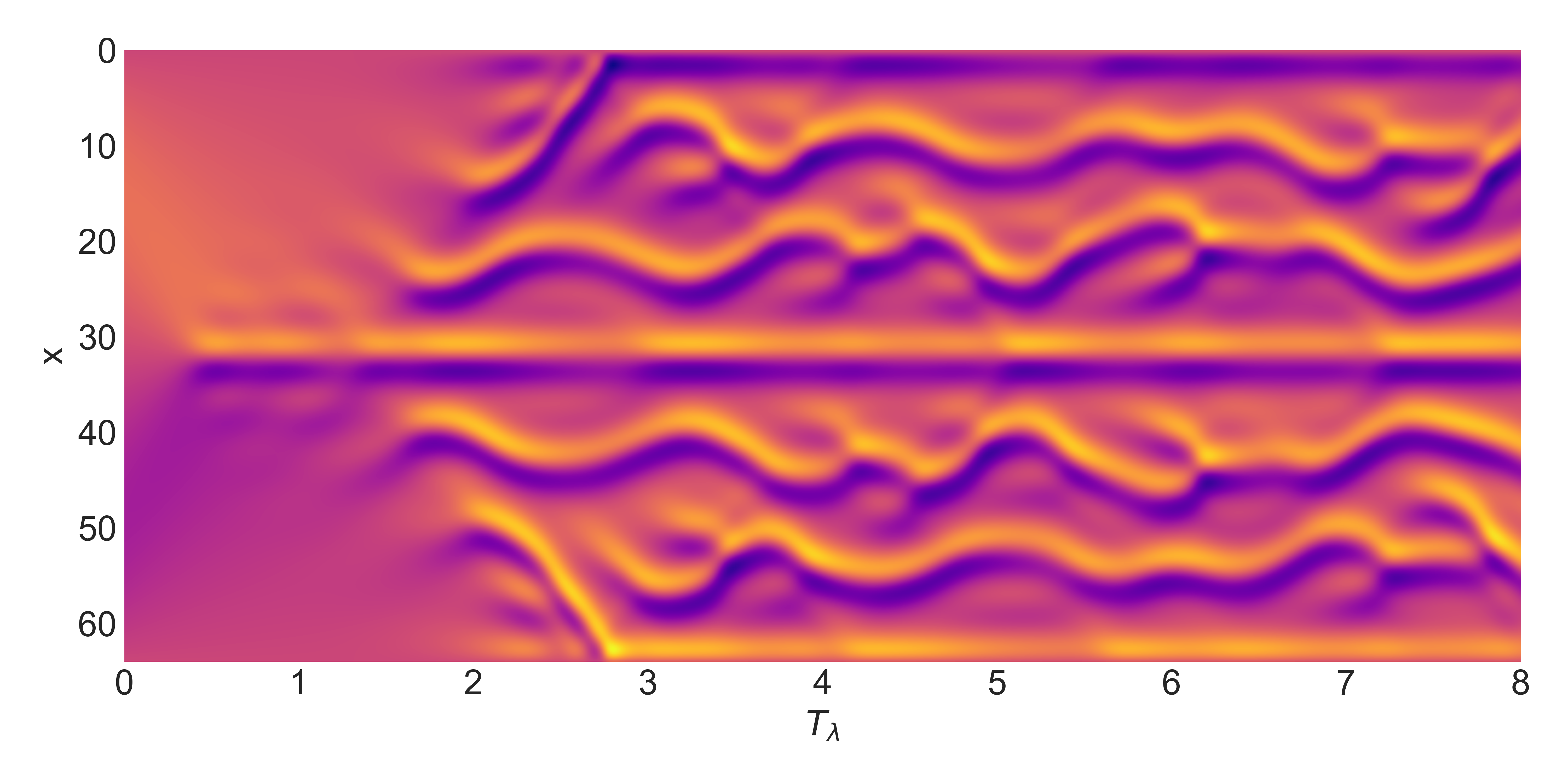}
    \caption{A numerical trajectory of the KS equation for 8 Lyapunov time with the initial condition $\vu(x, 0) = \sin(2\pi x / L)$. \label{fig:ks_imshow}}
\end{figure}

In one spatial dimension, the Kuramoto-Sivashinsky (KS) equation is given by
\begin{equation}\label{eq:ks}
    \vu_t = -\vu_{xx} - \vu_{xxxx} - \vu \vu_x,
\end{equation}
and is posed with periodic boundary condition $\vu(0, t) = \vu(L, t)$, for some length $L$.
This equation is widely studied as an archetypal example of a chaotic PDE and is considered one of the simplest such PDEs. It is also a useful surrogate for many fluid-dynamics applications, since it exhibits a wide range of complex dynamics including spatio-temporal chaos \cite{kuramoto1978diffusion}. The KS equation combines the linear antidiffusion and hyperdiffusion terms with a nonlinear advective term. Because of the structure of the linear terms, high frequency modes are stiffly damped by the hyperdiffusive part, while low frequency modes are \emph{excited} by the antidiffusive part. Although the linear part of this equation is apparently unstable, the nonlinear term stabilizes the equation.
The KS equation also features a non-trivial neutral manifold, corresponding with traveling waves. A typical trajectory for the KS equation over a chaotic time-scale is shown in Figure \ref{fig:ks_imshow}. The maximal LE for the KS equation depends on the length-scale $L$, which we fix at $L = 64$ for our experiments. With this parameter, we observe a maximal LE of approximately $0.1$, so one Lyapunov time $T_\lambda$ for the KS equation with these parameters is around 23 real time units. Being a stiff, nonlinear, and chaotic PDE with a non-trivial neutral stable manifold, the KS equation is a challenging problem to solve with PinT methods.

In the following experiments we discretize the KS equation using 4th order finite differencing in space and the fully implicit, two stage, 2nd order Lobatto IIIC method in time. We chose the Lobatto IIIC method since it is stiffly accurate and fairly simple to implement. The nonlinear implicit equation is solved using Newton's method, while the linear part of the Newton iteration is solved using UMFPACK from SuiteSparse \cite{umfpack}. When using MGRIT with naive time-coarsening and  rediscretization, we use the same Lobatto IIIC method with a larger time-step size. Since the two-stage Lobatto methods are actually a three parameter family of methods, including Lobatto IIIA (Crank-Nicolson), Lobatto IIIB, and Lobatto IIIC$^*$ (explicit trapezoid) which all use function evaluations calculated at the end-points, we can combine all four of these methods into a $\theta$ method. 
Solving algebraic equations analogous to those in Section \ref{sub:theta} allows us to find values of the three parameters, $\theta_A$, $\theta_B$, and $\theta_C$, which approximate the $m$-step stability function for Lobatto IIIC up to fourth order in $mh$. One degree of freedom is lost since we find the additional constraint that $\theta_A = \theta_B$. The resulting Runge-Kutta method is A-stable, and not stiffly accurate, and is thus a poor solver for the KS equation on coarse time-grids. So, we instead use one of the available degrees of freedom to enforce a stiff constraint on the stability function,
\begin{equation*}
    \lim_{z \to -\infty} \phi_\theta(z) = 0,
\end{equation*}
yielding a $\theta$ method of up to third order in $mh$ which is also stiffly accurate. 
In all cases where the $\theta$ method is used for the KS equation, this is the method used.  We use the initial condition $\vu(x, 0) = \sin(2\pi x / L)$, which is chosen since it is a smooth function which satisfies the periodic boundary condition and is easily generalized to any spatial grid resolution.

\subsubsection{Multilevel results}\label{sub:KS_scaling}
We use FCF-relaxation for the following scaling studies, since they improve convergence for MGRIT by relaxing parabolic modes of error (in the stable manifold) more effectively than F-relaxation~\cite{Do2016} and generally lead to faster time to solution in our experiments. Unless stated otherwise, a coarsening factor of $m=4$ is used. However in certain cases, we found that parallel efficiency could be improved by coarsening by a bigger factor (here $m = 16$) between the fine-grid and the first coarse-grid, while a coarsening factor of 4 was used for coarser grids, with little to no degradation in convergence, similar to~\cite{MGRIT14}. Further, as demonstrated in Section \ref{sub:lorenz_ml}, two-level convergence can be fast enough that the $\Delta$ correction is not needed, however, the $\Delta$ correction greatly improves the convergence for the multilevel method. Thus in practice, we defer the $\Delta$ correction to a coarser grid, while the first coarse-grid uses standard FAS MGRIT. This works well because the $\Delta$ correction is able to improve the accuracy of the coarsest grids where it is needed most, while significant work is saved by not computing $\Delta$ corrections on the finer grids where they are not as useful.  

In the following experiments, the $\Delta$ correction has been deferred to the second coarse-grid whenever more than two levels are used. An absolute residual stopping tolerance of $10^{-8}$ is used to terminate the MGRIT iterations.  Unless stated otherwise, the LV estimates are not propagated on the coarse-grid, as they are sufficiently resolved by fine-grid FCF-relaxation in most cases, as discussed in Section~\ref{sub:lowrank}.  MGRIT is an iterative method, and thus requires an initial guess for the solution across the whole time-domain. For the following experiments, the initial guess is produced by solving the coarsest grid with sequential time-stepping, then interpolating the solution to the fine-grid. Therefore, the quality of the initial guess is dependent on the accuracy of the coarsest grid.

\begin{figure}[h]
  \centering
  \includegraphics[width=0.56\textwidth]{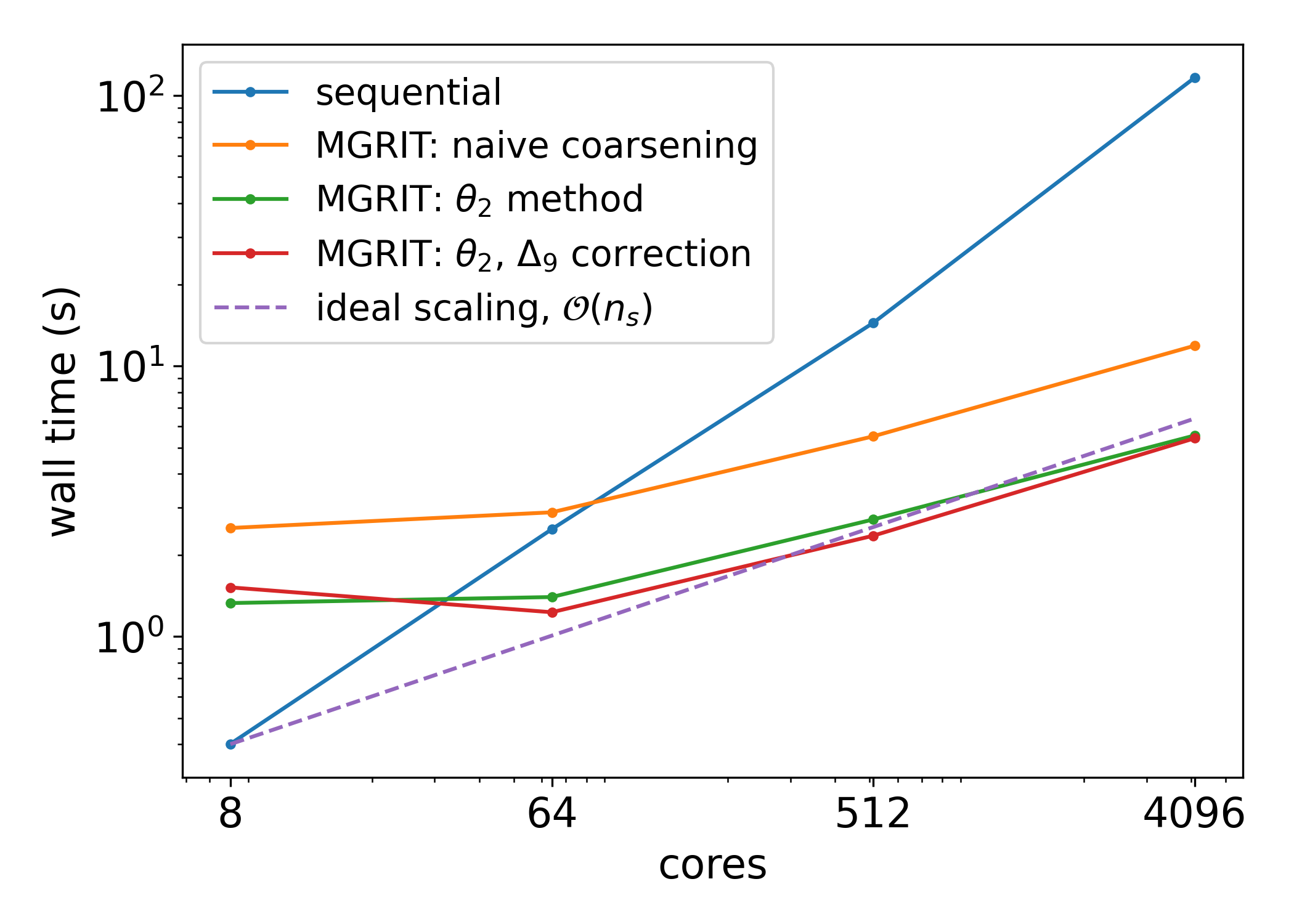}
  \caption{Weak scaling study for KS equation. \label{fig:ks_weak_refinement}}
\end{figure}

First, we demonstrate weak scaling for the KS equation.  Figure \ref{fig:ks_weak_refinement} plots the wall time to solution for increasing problem sizes, comparing sequential time-stepping against naive MGRIT, the $\theta$ method, and the $\theta$ method with rank 9 $\Delta$ correction. The KS equation is solved for 4 Lyapunov time (92.1 time-units), where the number of points in space are doubled as the points in time are quadrupled to maintain a fixed $h/h_x^2$ ratio and $h_x$ is the mesh size in space. The problem size per processor remains fixed at $256 \times 16$ points in space and time. The coarsest grid size is chosen to be as small as each algorithm stably supports for all problem sizes, remaining fixed at 128 time-points for naive MGRIT and MGRIT with the $\theta$ method, and 32 time-points for the $\theta$ method with rank 9 $\Delta$ correction. 
The ideal scaling plot is increasing because our spatial solver in Newton's method is UMFPACK, which has no parallelism in the spatial dimension and scales here for this 1D problem like $\mathcal{O}(n_x)$.
By choosing a suitable and stable coarsest grid size for naive MGRIT, naive MGRIT is able to scale well for a wide variety of problem sizes, albeit with a higher iteration count than for the other solvers.  This is reflected in the longer run-times for naive MGRIT. 
We believe that all algorithms benefit from the fact that the size of the unstable manifold is a property of the continuous equation, rather than being mesh dependent.
The $\theta$ method with rank 9 $\Delta$ correction achieves a max speedup of 21.5$\times$ over sequential time-stepping for the largest problem size, $(512 \times 32768)$.

\begin{figure}[h]
  \centering
  \includegraphics[width=0.56\textwidth]{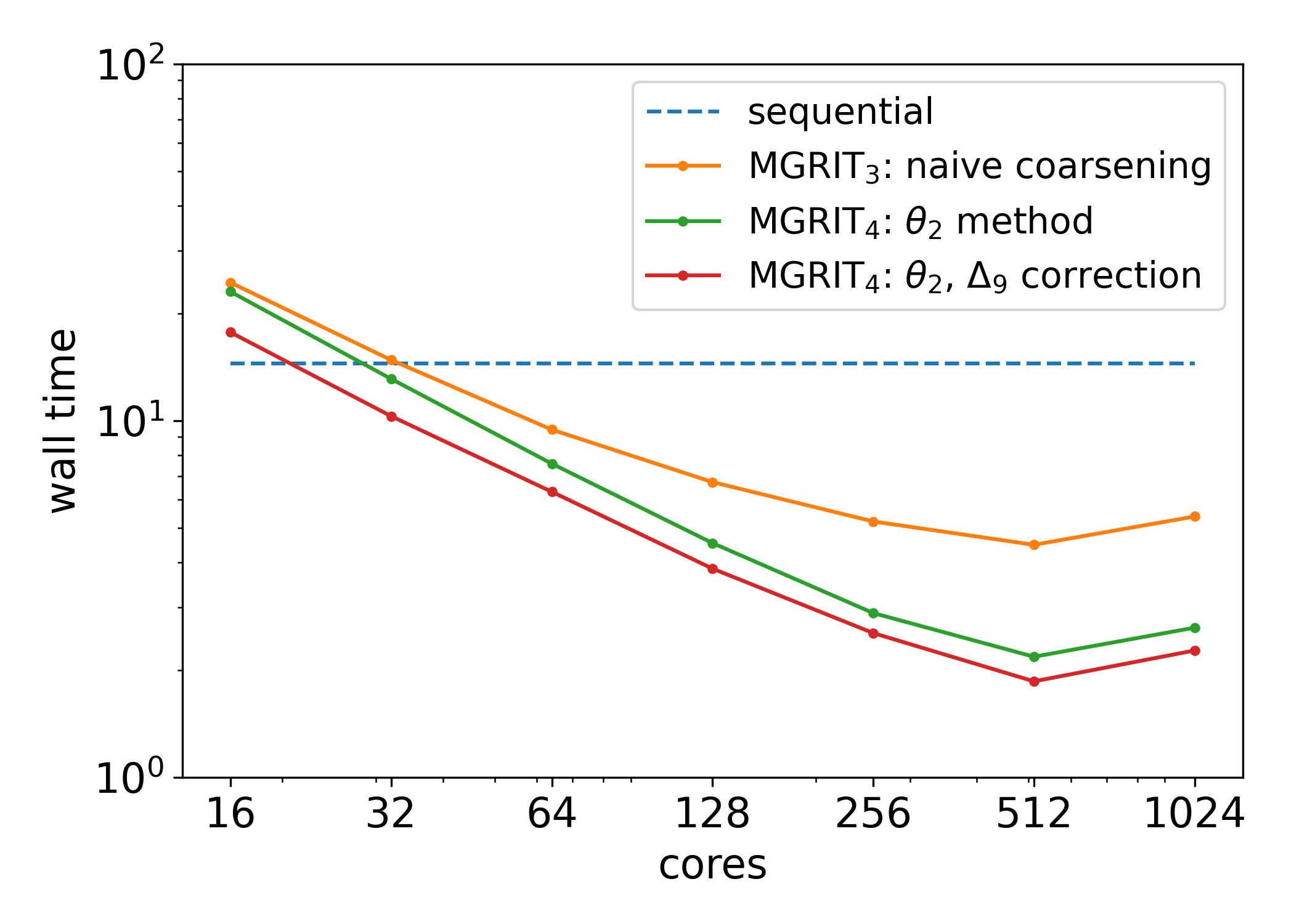}
  \caption{Strong scaling study for the KS equation with $T_f = 4T_\lambda$ and $(n_x, n_t) = (256, 16384)$, demonstrating a maximum speedup of 6.4$\times$ over sequential time-stepping. \label{fig:ks_strong}}
\end{figure}

\begin{figure}
    \centering
    \includegraphics[width=0.56\textwidth]{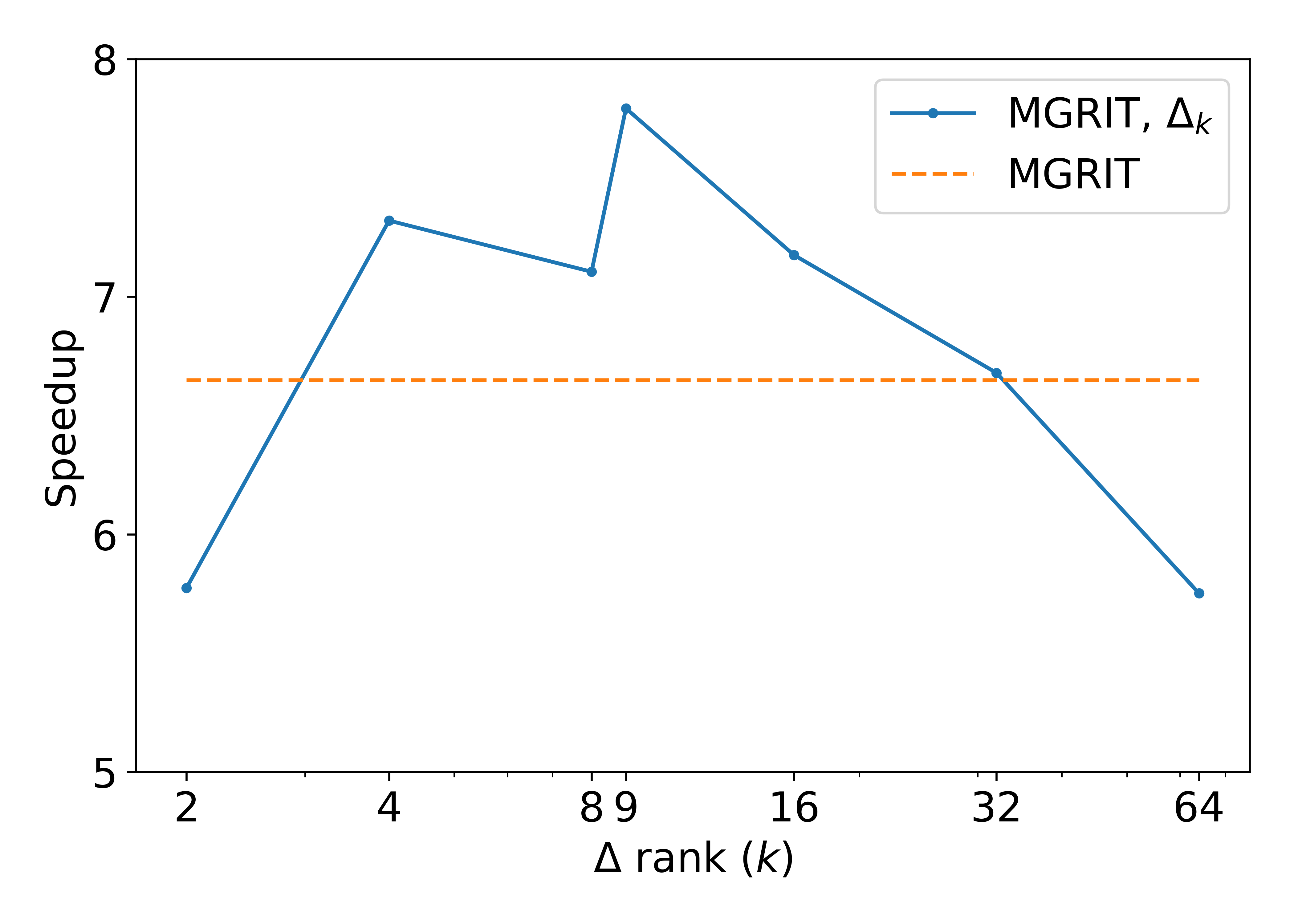}
    \caption{Maximum speedup over sequential time-stepping with 512 cores for the KS equation with varying ranks of low-rank $\Delta$ correction (solid line). As compared to MGRIT with the $\theta$ method coarse-grid (dashed line).  The best speedup is achieved when the rank is equal to the dimension of the unstable manifold (9). \label{fig:rank_dependence}}
    
\end{figure}

Next we explore strong scaling for the KS equation. Figure \ref{fig:ks_strong} plots the wall-clock time to solution for naive MGRIT, the $\theta$ method, and $\theta$ method with varying ranks of $\Delta$ correction, as compared to sequential time-stepping.  The KS equation is solved for a time-domain of 4 Lyapunov time (92.1 time-units) with 256 points in space and 8192 points in time. This corresponds with the second largest problem size from Figure~\ref{fig:ks_weak_refinement}. Naive MGRIT is limited to a 3-level method, since the method becomes unstable for 4-levels, and thus the sequential cost of solving the coarsest grid leads to quick run-time stagnation. The $\theta$ method, however, provides a stabilized coarse-grid, thus enabling a 4-level method, which scales better than naive MGRIT. The low-rank $\Delta$ correction reduces the number of iterations required for convergence with only a rank 9 correction, thus making the method slightly faster, and we observe a maximum speedup of 6.4$\times$. However, increasing the rank of the $\Delta$ correction beyond 9 does not further improve convergence, and only increases the cost of each iteration. Thus rank 9 $\Delta$ correction appears to be optimal for this problem, which is expected, as experiments measuring the full Lyapunov spectrum indicate that the unstable manifold is roughly 9-dimensional. The dependence of the maximum speedup on the rank of $\Delta$ correction is illustrated in Figure~\ref{fig:rank_dependence}. 

\begin{figure}[h]
    \centering
    \includegraphics[width=0.56\textwidth]{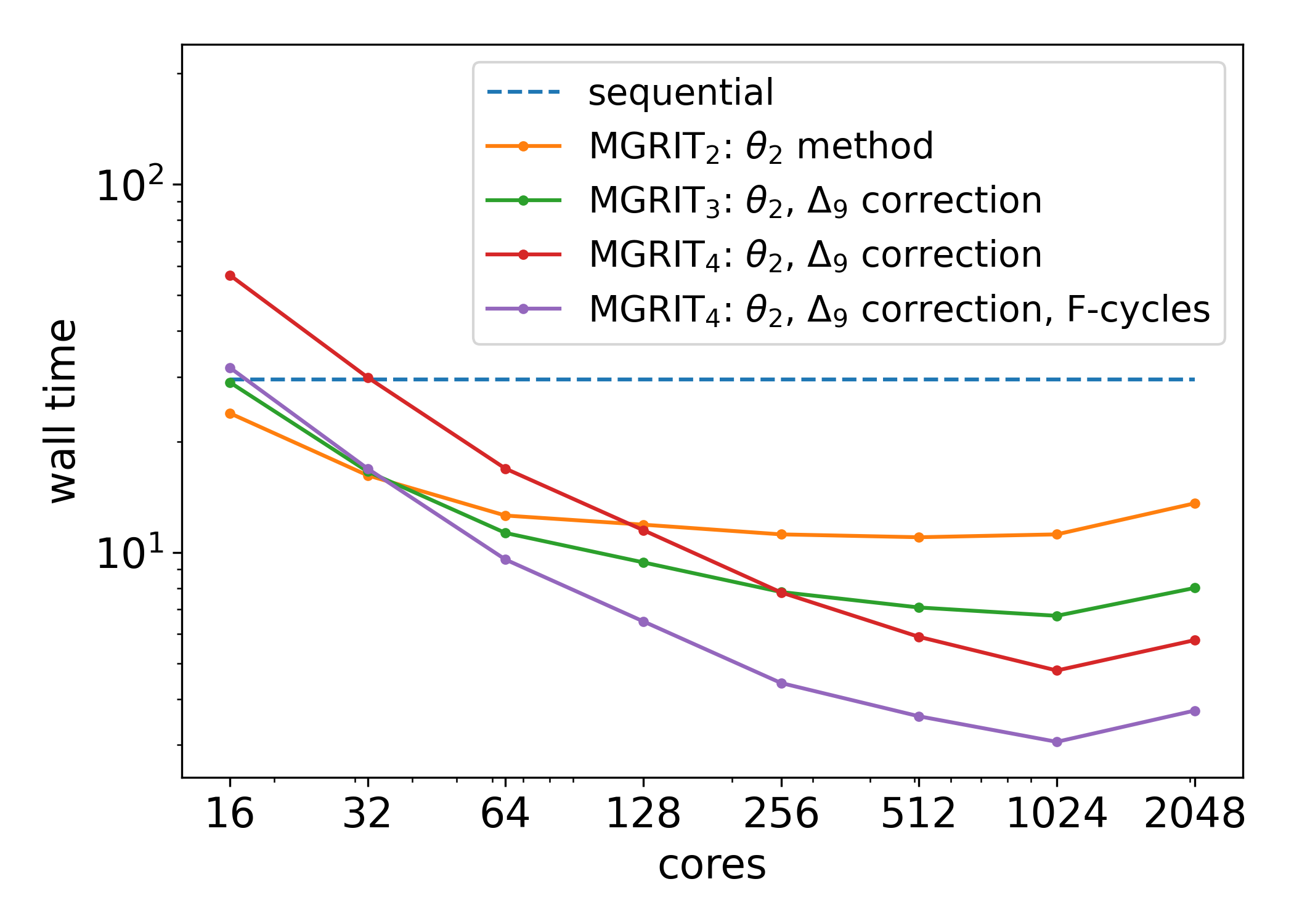}
    \caption{Strong scaling study for the KS equation with $T_f = 8T_\lambda$ and $(n_x, n_t) = (256, 16384)$, demonstrating a maximum speedup of 9.6$\times$ when using a four-level F-cycle with the $\theta$ method and rank 9 $\Delta$ correction.  \label{fig:ks_strong_tf8}}
\end{figure}

Last, we present a strong scaling study for an even longer time-domain. Figure~\ref{fig:ks_strong_tf8} demonstrates strong scaling for the KS equation over 8 Lyapunov time (184.2 time-units). As before, the KS equation is solved with 256 points in space, but the time-grid is increased to 16384 points in time, to match the fine-grid time-step size from the previous problem. For this problem, naive MGRIT is limited to a two-level method, with a coarse-grid of 1024 time-points, even with the $\theta$ method. This is because on this time-scale, MGRIT begins to stall when a third level is added. However, a rank 9 $\Delta$ correction is able to remedy this, allowing rapid convergence for both the three and four level methods. For this experiment, accuracy for the LV estimates requires that the coarse-grid propagation of the LVs be switched on (i.e., the Gram-Schmidt step in Algorithm \ref{alg:lrDelta}), otherwise the method stalls. Here, $\Delta$ correction allows for a dramatic improvement to the scaling of the method, and we observe a maximum speedup of 6.1$\times$. If we use F-cycles, which can improve convergence at the cost of some parallel efficiency \cite{MGRIT14}, the four-level method can be made even faster, achieving a maximum speedup of 9.6$\times$ for this problem.

\section{Conclusions}%
\label{sec:conclusions}
Solving chaotic dynamical systems with PinT is inherently difficult due to the exponentially increasing condition number with $T_f$. 
However, given that future speedups will require increased concurrency in codes, PinT is nonetheless needed for chaotic systems.

In this work, we dramatically improved the convergence of MGRIT for the considered chaotic problems by introducing a modified time-coarsening scheme (called $\theta$ method) and modified FAS coarse-grid (called $\Delta$ correction).  The $\Delta$ correction method is equivalent to applying the MNM method, originally designed for nonlinear elliptic equations, to the time-dimension. We further introduced a novel low rank (low memory) approximation to $\Delta$ correction which is needed for problems with large spatial dimensions.  The results are promising, where for some cases, we observed nearly flat iteration counts for long time-domains.  For the KS equation, we observed meaningful PinT speedup for a chaotic problem (up to nearly 10$\times$ speedup over sequential time-stepping), which is the first such speedup to our knowledge.

Based on our results, we expect that replacing the spatial Newton solver with a more robust, parallel nonlinear solver (such as classical FAS in space or the MNM algorithm) would allow us to coarsen even further in time. This would improve the strong scaling results presented in \ref{sub:KS_scaling} to better match the optimal scaling demonstrated for linear parabolic equations in \cite{MGRIT14}. Also, using an optimally scaling parallel solver for the spatial dimension, combined with the low rank $\Delta$ correction, would likely result in an algorithm with mesh-independent weak scaling. This is future work.

Finally, although the solution to the IVP is the focus of this paper, if the estimated Lyapunov vectors generated by the low-rank $\Delta$ correction algorithm are high enough quality, the presented algorithm has the potential to \emph{greatly accelerate the Lyapunov analysis of chaotic systems} by simultaneously solving for the state solution and the LVs, parallel in time. This would represent a large improvement over current algorithms which are sequential in time and prohibitively expensive for many applications.

\bibliographystyle{siamplain}
\bibliography{references_no_url}
\end{document}